\documentclass[12pt]{amsart}
\usepackage{graphicx,amssymb,amsfonts,epsfig,amsthm,a4,amsmath,url}
\usepackage{mathrsfs}
\usepackage{fancyhdr}
\usepackage[lined,boxed,commentsnumbered]{algorithm2e}
\usepackage{epstopdf}
\usepackage[font={small,it}]{caption}
\usepackage{subcaption}
\usepackage[section]{placeins}
\input amssym.def
\input amssym
\usepackage{hyperref}
\hypersetup{
    pdftitle={Dynamical sampling with additive random noise},    
    pdfauthor={Akram Aldroubi}{Ilya Krishtal},     
    pdfsubject={Krylov subspace methods in dynamical sampling},   
    pdfcreator={Ilya Krishtal},   
    pdfproducer={Ilya Krishtal}, 
    pdfkeywords={.}{.}{.}, 
     colorlinks,%
    citecolor=blue,%
    filecolor=black,%
    linkcolor=magenta,%
    urlcolor=black
}

\chardef\bslash=`\\ 

\makeatletter
\def\verbatim{\interlinepenalty\@M \@verbatim
   \leftskip\@totalleftmargin\advance\leftskip2pc
   \frenchspacing\@vobeyspaces \@xverbatim}
\makeatother
\newtheorem{thm}{Theorem}[section]
\newtheorem{cor}[thm]{Corollary}
\newtheorem{lem}[thm]{Lemma}
\newtheorem{prop}[thm]{Proposition}

\theoremstyle{defn}
\newtheorem{defn}{defn}[section]

\theoremstyle{rem}
\newtheorem{rem}{Remark}[section]
\newtheorem{exmp}{Example}[section]

\numberwithin{equation}{section}

\newcommand{\begeq}{\begin {equation}}
\newcommand{\eq}{\end{equation}}
\newcommand{\bs}{\begin {split}}
\newcommand{\es}{\end{split}}
\newcommand{\bp}{\begin {prop}}
\newcommand{\ep}{\end {prop}}
\newcommand{\bt}{\begin {thm}}
\newcommand{\et}{\end {thm}}
\newcommand{\bc}{\begin {cor}}
\newcommand{\ec}{\end {cor}}
\newcommand{\bl}{\begin {lem}}
\newcommand{\el}{\end {lem}}
\newcommand{\bpf}{\begin {proof}}
\newcommand{\epf}{\end {proof}}
\newcommand{\bi}{\begin {itemize}}
\newcommand{\ei}{\end {itemize}}
\newcommand{\ben}{\begin {enumerate}}
\newcommand{\een}{\end {enumerate}}
\newcommand{\brem}{\begin {rem}}
\newcommand{\erem}{\end {rem}}
\newcommand{\bex}{\begin {exmp}}
\newcommand{\eex}{\end {exmp}}
\newcommand{\bd}{\begin {defn}}
\newcommand{\ed}{\end {defn}}

\newcommand{\ZZ}{{\mathbb Z}}

\newcommand{\RR}{{\mathbb R}}
\newcommand{\CC}{{\mathbb C}}

\newcommand{\vecx}{\bar{\mathbf{f} }}
\newcommand{\vecy}{\bar{\mathbf{y} }}

\allowdisplaybreaks

\begin{document}
	%
	\title{Dynamical sampling with additive random noise}

	\author[A. Aldroubi et.~al.]{Akram Aldroubi,
			Longxiu Huang, Ilya Krishtal,\\
			 Akos Ledeczi, Roy R. Lederman, Peter Volgyesi}
			
%

\date{\today }



	\maketitle

\begin{abstract}
	Dynamical sampling deals with   signals  that evolve in time under the action of a linear operator. The  purpose of the present paper is to analyze the performance of the basic dynamical sampling algorithms in the finite dimensional case and study the impact of additive noise.   The  algorithms are implemented and tested on  synthetic and real data sets, and denoising techniques are integrated to mitigate the effect of the noise.   We also develop theoretical and numerical results that validate the algorithm for recovering  the driving operators, which are defined via a real symmetric convolution.
	\vspace{3mm}
\end{abstract}

%
%

\section{Introduction}

Dynamical sampling is a framework for processing signals that evolve in time under the action of a linear operator. In dynamical sampling, one seeks to exploit the association between
the signals received at various time levels to enhance the classical sampling and reconstruction techniques or propose novel sampling and reconstruction algorithms. Since the original work on dynamical sampling \cite{ADK13}, a number of subsequent studies have been devoted to various aspects of the theory and  applications (see, for example,  \cite{APT15,AT14, ACCMP17, ACMT17, ADK15, AK16,  AKT17, AKW15,CMPP17, CM17, D14, GRUV15, JD15,JD17,  P17, Phi17, Tan15, ZLL17_2, ZLL17}).

The present study addresses certain numerical and theoretical aspects of the two main problems of dynamical sampling in the finite dimensional setting. The first problem is to recover a signal $f$ that evolves in time  under the action of a known operator $A$ \cite{ACMT17, ADK13}. The second problem is to recover the driving operator $A$ in the case when it is unknown or only partially known \cite{AK16}.  The main contributions of this study are as follows: (a) performance evaluation  of basic algorithms developed within the dynamical sampling framework; (b) analysis of the impact of additive noise and the effectiveness of denoising techniques,  when processing real and synthetic data sets. A preliminary version of this study has been documented in \cite{AHKL17}.

As we mentioned above, the first  problem of dynamical sampling is concerned with the recovery of a signal $f$  that evolves in time  under the action of a known operator $A$.
More precisely, we consider a signal $f \in \mathbb{C}^d$ and a bounded  linear operator $A$ on $\mathbb{C}^d$ which we identify with its matrix in the standard basis. At  time level $n\in \mathbb N$,  the signal becomes
\begin{equation}\label{equation1}
f_n=A^nf. 
\end{equation} 
We let $\Omega\subset\{1,\ldots,d\}$ denote a set of ``spatial" locations. 
The noiseless dynamical samples  are then 
\begin{equation}
\label {DynSamp}
\left\{f_n(j):j\in\Omega, 0\le n\le L \right\}.
\end{equation}  
In  \cite{ACMT17}, necessary and sufficient conditions for recovering $f \in \mathbb{C}^d$ have been derived in terms of $A$, $\Omega$, and $L$.

In the noisy case, we consider the corrupted dynamical samples of the form 
\begin {equation}
\label {NoisyDySamp}
\left\{f_n(j)+\eta_n(j),\ j\in \Omega,\ 0\le n\le L\right\},
\end{equation}
where $\eta_n$, $n\ge0$   are independent identically distributed (i.i.d.) $d$-dimensional random variables  with zero mean   and covariance matrix  $\sigma^2I$, and  $\eta_n({j})$ denotes the $j$-th component of $\eta_n$.

Using  the  $d\times d$ diagonal sub-sampling matrix $S_{\Omega}$, defined by 
\begin {equation}
\label {SamOp}
(S_{\Omega})_{jj}=\begin{cases} 
1 &  j\in\Omega \\
0 & \text{otherwise},
\end{cases}
\end {equation}
the noisy  data sampled at time level $n$ in \eqref {NoisyDySamp} can be described by the vectors $\tilde{y}_n$ given by
\begin {equation}
\tilde{y}_n = S_{\Omega}(f_n+\eta_n)
\end{equation}

The signal $f$ can be approximately recovered from the noisy measurments $\tilde{y}_n$ by solving the least squares minimization problem
\begin{eqnarray}\label{equation 1.2}
f^\sharp_L=\arg\min_{g}\sum_{n=0}^{L}\left\|S_{\Omega}(A^ng)-\tilde{y}_n\right\|^2_2.\label{equation2}
\end{eqnarray}
In this study, an iterative algorithm for solving problem \eqref{equation 1.2} is investigated. In addition, the mean squared error (MSE) $E(\|\epsilon_L\|^2)$ is estimated with $\epsilon_L=f^\sharp_L-f$ and  the behavior of the MSE is analyzed  as $L\rightarrow\infty$ for an unbiased linear estimator. 

The second  problem of dynamical sampling deals with the case when the evolution operator $A$ is unknown (or only partially known). In \cite{AK16}, an algorithm has been proposed for finding the spectrum of $A$ from the dynamical samples. The present paper delves deeper into this algorithm from both theoretical and numerical perspectives. From the theoretical perspective, an alternative proof is given for the fact that the algorithm in \cite{AK16} can  (almost surely)  recover the spectrum of $A$ from dynamical samples and also recover the operator $A$ itself, in the case when it is known that $A$ is given by a circular convolution $Af=a*f$ with some real symmetric filter $a$ in $\RR^d$. From a numerical  point of view, this analytical result lays the theoretical foundation and paves the way toward recovering the operator $A$ and the signal from real data collected from physical processes such as the heat diffusion.   The nature of the spectrum recovery algorithm also motivates an integration of  Cadzow-like denoising techniques \cite{C80, G10c}, which can be applied to both synthetic  and real data. 

\subsection {Contribution and Organization} 
In Section \ref{prelim}, we summarize the notation that is used throughout the paper and present the algorithms for signal and filter recovery that work ideally in the noiseless case. To recover the signal, we borrow a least squares updating technique  from \cite{B84LS} and  tailor it for dynamical sampling. To recover the driving operator (in the case of a convolution), we review the  algorithm from \cite{AK16} and provide its new  derivation, which is more straightforward than the general proof in \cite{AK16}.
In Section \ref{noise_red}, the Cadzow denoising method  is sketched for a special case of uniform sub-sampling; it is validated to be numerically efficient in the context of dynamical sampling   in Section \ref{test}. Section \ref{errant} is dedicated to the error analysis of the least squares solutions for finding the original signal in the presence of additive white noise.
It shows the relation between the MSE of the solution and the number of time levels considered. In Section \ref{test}, we outline the outcomes of the extensive tests performed for the algorithms discussed in Sections \ref{prelim} and \ref{noise_red}. More precisely, Section \ref{sim1} demonstrates  the consistency of the theory for the MSE of the least squares solutions on synthetic data. Section \ref{cd} illustrates the effect of Cadzow denoising method on signal and filter recovery in the case of synthetic data. Finally, in Section \ref{rd}, the recovery algorithms and denoising techniques are integrated together to process real data collected from cooling processes.

\section {Notation and Preliminaries}\label{prelim}
\subsection {Notation}
Let $\mathbb{Z}$ be the set of all integers and $\ZZ_d$ be the cyclic group of order $d$. By $\CC^d$ and $\CC^{m\times d}$ we denote the linear space of all column vectors with $d$ complex components and the space of complex matrices of dimension $m\times d$, respectively. Given a matrix $A\in\mathbb{C}^{m\times d}$, $A_{ij}$ stands for the entry of the $i$-th row and $j$-th column of $A$, $A^*$ represents the conjugate transpose of $A$, and the 2-norm of $A$ is defined by
$$\|A\|=\sup_{f\in\mathbb{C}^d, \|f\|_2=1}\|Af\|_2,$$
where $\|f\|_2=\sqrt{\sum_{i=1}^{d}|f(i)|^2}$ and $f(i)$ refers  to $i$-th component of a vector $f \in \CC^d$. 

For a random variable $x$ that is distributed normally with mean $\mu$  and variance $\sigma ^{2}$,  we may write 
$x\sim N(\mu,\sigma^2).$
\subsection{A general  least squares updating technique for signal recovery} 
We borrow  from \cite {GVL13} the following updating technique  for adjusting a least squares solution when new equations are added.
Consider the following least squares problem
\begin{equation} 
\label {GLS}
f^\sharp_L=\arg\min_{g\in\mathbb{C}^d}\sum_{i=1}^{L}\left\|A_ig-b_i\right\|_2^2,
\end{equation}
where $A_i\in\mathbb{C}^{m_i\times d}$, and $\text {rank} (A_1)=d$ (i.e., $A_1$ has  full column rank).   

We take the case of $L=2$ as an example to explain the updating technique.
Consider the QR decomposition  $A_1=Q_1R_1$, where $Q_1$ is an $m_1\times d$  matrix satisfying $Q_1^*Q_1=I$ and $R_1$ is a $d\times d$ triangular matrix. Then $$f^\sharp_1=\arg\min_{g}\left\|A_1g-b_1\right\|_2^2=\arg\min_{g}\left\|R_1g-Q_1^*b_1\right\|_2^2.$$ Let $\tilde{b}_1=Q_1^*b_1$. Since $A_1$ has full rank,  we have  $f^\sharp_1=R_1^{-1}\tilde{b}_1.$
Suppose that new information is added, then the least squares problem and its solution needs to be updated, i.e., $f^\sharp_2=\arg\min\limits_{g\in\mathbb{C}^d}\sum_{i=1}^{2}\left\|A_ig-b_i\right\|_2^2.$   

To solve the new least squares problem, we note that
\begin{eqnarray*}
\arg\min_{g\in\mathbb{C}^d}\sum_{i=1}^{2}\left\|A_ig-b_i\right\|_2^2&=&\arg\min_{g\in\mathbb{C}^d}\left\|\left(\begin{array}{c}
A_1\\
A_2
\end{array}\right)g-\left(\begin{array}{c}
b_1\\
b_2
\end{array} \right)
\right\|_2^2\\
&=&\arg\min_{g\in\mathbb{C}^d}\left\|\left(\begin{array}{cc}
Q_1& 0\\
0& I
\end{array}\right)\left(\begin{array}{c}R_1\\
A_2
\end{array}\right)g-\left(\begin{array}{c}
b_1\\
b_2
\end{array}\right)\right\|_2^2\\
&=&\arg\min_{g\in\mathbb{C}^d}\left\|\left(\begin{array}{c}
R_1\\
A_2
\end{array}\right)g-\left(\begin{array}{c}
Q_1^*b_1\\b_2
\end{array}\right)\right\|_2^2\\
&=&\arg\min_{g\in\mathbb{C}^d}\left\|\left(\begin{array}{c}
R_1\\
A_2
\end{array}\right)g-\left(\begin{array}{c}
\tilde{b}_1\\b_2
\end{array}\right)\right\|_2^2. 	
\end{eqnarray*}
Therefore, the problem reduces to finding $$f^\sharp_2=\arg\min_{g\in\mathbb{C}^d}\left\|\left(\begin{array}{c}
R_1\\
A_2
\end{array}\right)g-\left(\begin{array}{c}
\tilde{b}_1\\b_2
\end{array}\right)\right\|_2^2. $$
One further needs to calculate the  QR decomposition  $$\left(\begin{array}{c}
R_1\\
A_2\end{array}\right)=Q_2R_2,$$ where $Q_2$ is a unitary matrix and $R_2$ is a $d\times d$ triangular matrix. Denote $$\tilde{b}_2=Q_2^*\left(\begin{array}{c}
\tilde{b}_1\\
b_2
\end{array}\right).$$ It follows that $f^\sharp_2=R_2^{-1}\tilde{b}_2.$

The same process can be applied to the case $L\geq 3$ which leads to the iterated updating algorithm that is summarized in Algorithm \ref{algorithm iterative}.

This algorithm demonstrates that the recovery problem in dynamical sampling can be solved in a streaming setup,
where the solution is updated as new measurements are collected over time, 
\begin{enumerate}
\item[1.] without storing all the previous samples ($b_j$) or explicitly rewriting all the matrices ($A_j$) for all $j<i$ at the $i$th step,
\item[2.] and taking advantage of quantities that are stored from previous iterations to avoid the naive computation involving all the previous samples and matrices. 
\end{enumerate}

Observe that in the dynamical sampling framework we have $A_i=S_{\Omega}A^{i-1}$. Assume that at step $i$, the  QR decomposition for 
$$\mathcal{A}_{i}=\left(\begin{array}{c}
S_{\Omega}I\\
S_{\Omega}A\\
\vdots\\
S_{\Omega}A^{i-1}
\end{array}\right) $$
is $$\mathcal{A}_i=QR.$$
At   step $i+1$, 
$\mathcal{A}_{i+1}$ can thus be written in the convenient form $$\left(\begin{array}{c}
S_{\Omega}I\\
QRA
\end{array}\right)=\left(\begin{array}{cc}
I &0\\
0&Q
\end{array}\right)\left(\begin{array}{c}
S_{\Omega}I\\
RA
\end{array}\right).$$

\IncMargin{1em}
\begin{algorithm}[H]\label{algorithm iterative}
\SetKwInOut{Goal}{Goal}

\Goal{Recover the original signal by processing time series data.}
{\text {\bf Input} $A_1, b_1$}

{Set $A_1=Q_1R_1$,  the economic QR decomposition of $A_1$ with the assumption that $A_1$ has full column rank (see Remark \ref {ChoiceOmega})}.

{Set $\tilde{b}_1=Q_1^*b_1$.}

{
{Set  $f^\sharp_1=R^{-1}\tilde{b}_1$.}}

\For{$i = 2$ \KwTo $L$}
{ {\text {\bf Input} $A_i, b_i$}

Compute the QR decomposition for $\left(\begin{array}{c}
R_{i-1}\\
A_i
\end{array}\right) =Q_iR_i$ using the Householder transformation \cite {GVL13}.\\
Set
$\tilde{b}_i=Q_i^*
\left(\begin{array}{c}
\tilde{b}_{i-1}\\
b_i
\end{array}\right)$.

{
	{Set $f^\sharp_i=R^{-1}_i\tilde{b}_i$.}}

}
\text{\bf Output} {$f^\sharp_L$}
\caption{Pseudo-code of the iterated updating algorithm.}
\end{algorithm}

\subsection {Filter recovery for the special case of convolution operators and uniform subsampling}\label{frec} 
In this section, we recall from \cite{AK16} an algorithm  for recovering an unknown driving operator $A$ that is defined  via a convolution with a   real symmetric filter i.e., $A$ is a circulant matrix corresponding to a convolution with $a$: $Af=a\ast f$), and where the spatial sampling is uniform at every time-instant $n$. 
We also provide a new, direct proof of validity for the filter recovery algorithm for this case.  
Specifically, we consider samples of $A^\ell f=a^\ell\ast f$  at $m\ZZ_d$ where $m\ge 2$, and  $a^
\ell=a\ast\cdots\ast a$ is the $\ell$ times convolution of the filter $a$.
We also assume that the Fourier transform $\hat a$ of the filter $a$   is real symmetric, and strictly decreasing on $[0,\frac {d-1} 2]$.
We will use the notation $S_mf_n$ to describe this uniform subsampling. In particular, for a vector $z \in \ell^2{(\ZZ_d)}$, $S_mz$ belongs to $\ell^2(\ZZ_{J})$, and $S_mz(j)=z(mj)$ for $j=1,\dots,J$, where throughout we will assume that $m$ is odd, and $d=Jm$ for some odd integer $J$. 

Let  
\begeq\label{data}
y_\ell = S_m(A^\ell f)=S_m(a^\ell\ast f), \; \ell\geq 0,
\eq
be 
the dynamical samples at time level $\ell$. 
By Poisson's summation formula, 
\begeq \label {DSformula} \widehat{(S_mz)}(j) = \frac{1}{m}\sum\limits_{n=0}^{m-1}\hat{z}(j + nJ), \quad 0\leq j\leq J-1, \ z \in \ell^2(\ZZ_d),
\eq 
An application of the Fourier transform to  \eqref {data} yields
\begeq \label {datarel}
\hat y_\ell(j)= \frac{1}{m}\sum\limits_{n=0}^{m-1}\hat{a}^\ell(j + nJ)\hat {f}(j+nJ), \quad 0\leq j\leq J-1.
\eq
For each fixed $j \in\ZZ_J$ and for some integer $L$ with $L\ge 2m-1$  ($L=2m-1$ is the minimum number of time levels that we need to recover the filter), we introduce the following notation: $$\vecy_\ell (j)=\big({\hat y_{\ell}(j),\hat y_{\ell+1}(j),\dots,  \hat y_{\ell+L}(j)}\big)^T,$$ 
$$\vecx (j)=\big({\hat f(j),\hat f(j+J),\dots,  \hat f(j+(m-1)J)}\big)^T,$$ and 
\begeq \label {fdNScond} \mathcal V_m(j)= \left(
\begin{array}{cccc} 1 & 1 & \hdots &1\\
\hat{a}(j) & \hat{a}(j+J) & \hdots & \hat{a}(j+(m-1)J)\\
\vdots & \vdots& \vdots & \vdots \\
\hat{a}^{L-1}(j) & \hat{a}^{L-1}(j+J) & \hdots & \hat{a}^{L-1}(j+(m-1)J)\end{array}\right),  
\eq
where $0\leq j\leq J-1. $
From \eqref {datarel},  it follows that 
\begeq \label {Fouriereq0}
\vecy_\ell (j) = \frac1{m} \mathcal V_m(j)D^\ell(j) \vecx (j), \text{ for } 0\leq j\leq J-1, ~\ell\geq 0,
\eq
where $D(j)$ is the diagonal matrix  $D(j)=\text{diag} \big(\hat{a}(j), \hat{a}(j+J) , \dots,\hat{a}(j+(m-1)J)\big)$. Let  $p_j(x)=c_0(j)+c_1(j)x+\dots+ c_{1}(j)x^{n_j-1}+x^{n_j}$ be the minimal polynomial that annihilates 
$D(j)$. The degree of $p_j$ is equal to the number of distinct diagonal values of $D(j)$. 
Since $L \geq 2m-1$, it follows from the assumptions on $\hat a$ ($\hat a$ is real symmetric, and strictly decreasing on $[0,\frac {d-1} 2]$)
that $\deg (p_j)=m$ for $j\ne 0$ and $\deg (p_0)=(m+1)/2$.  Moreover, the rectangular Vandermonde matrix $\mathcal V_m(j)$ has rank $r_j=m$ if $j\ne 0$, and  $r_0=(m+1)/2$ if $j=0$. Consequently, using \eqref {Fouriereq0}, we have that  for almost all $\hat f$,
\begeq \label {Fouriereq}
\vecy_{k+r_j}(j) + \sum_{\ell=0}^{r_j-1}c_\ell (j) \vecy_{k+\ell}(j) =0,  \quad 0\leq j\leq J-1,
\eq
where $c_\ell (j)$ are the coefficients of the polynomial $p_j$ and  $r_j=\deg p_j= 
\text {rank } \mathcal V_m(j) $. 
The above discussion leads to  the following Algorithm  \ref{spectrum recovery for convolution operators} for recovering the spectrum $\sigma(A)$. 

\IncMargin{1em}
\begin{algorithm}[H]\label{spectrum recovery for convolution operators}
\SetKwInOut{Goal}{Goal}

\Goal{Recover the spectrum $\sigma(A)$.} 
\BlankLine
Set $J=d/m$.
\\
\For{$j = 0$ \KwTo $J-1$}{
Find the minimal integer $r_j$ for which the system \eqref{Fouriereq} has a solution $c(j)$ and find the solution\;
set $p_j(\lambda)=\lambda^{r_j}+\sum_{\ell=0}^{r_j-1}c_\ell(j)\lambda^\ell$ and find the set $R(j)$ of all roots of $p_j$.
}
{Set $\sigma(A)=\bigcup_{j=0}^{J-1} R(j)$.}
\caption{A spectrum recovery algorithm for convolution operators.}
\end{algorithm} 
\begin{rem} The algorithm for spectrum recovery involves finding the roots of a set polynomials of degree $m$ or $\frac {m+1} 2$, where $m$ is the subsampling factor.  This problem becomes more and more difficult as $m$ becomes larger and larger. However, in applications, one could expect $m$ to be of moderate size ($m\le 5$). Moreover, if some of the spectral values are too close to each other, then finding the coefficients of the minimal polynomials becomes  unstable.
\end {rem}

\begin{rem}
The recovery of both the filter and the signal from the measurements points to certain relations to the problem of
Blind Deconvolution (see for example \cite{kundur1996blind}); typically, Blind Deconvolution does not involve the difficulty arising from the sub-sampling (the operator $S_m$), but it is restricted to one time measurement, and uses other assumptions on the signal and filter.
\end{rem}

\section {Cadzow Denoising Method} \label{noise_red}
\label{Cadzow denoising theory} In this section, we describe a Cadzow-like algorithm (see Algorithm \ref{Cadzow denoising algorithm}) \cite{C80, G10c} which can be effectively applied to  approximate   the dynamical samples  $y_n$ in \eqref {data} from   the noisy measurements $\tilde {y}_n= y_n+\eta_n$. 

Suppose data points $y_n$ in \eqref{data} are  such that  $m$ is an odd integer and $A$ is a symmetric
circulant matrix generated by a real symmetric filter $a$, i.e., the Fourier transform $\hat a$ of the filter $a$ is real symmetric. In addition, we also assume that   $\hat a$   is  monotonic on $[0, \frac{d-1}2]$. 
Let $L$ be the number of time levels as in \eqref {fdNScond}. In particular, it is necessary that $L\ge 2m-1$. Without loss of generailty we assume that $L$ is even. From \eqref {datarel}, \eqref {Fouriereq0}, \eqref {Fouriereq}  in Section \ref{frec} 
(see also \cite{ADK13, AK16}), it follows that the Hankel  matrix
\begin{equation}\label{Hj}
H(j)=\left(\begin{array}{cccc}
\hat y_0(j)&\hat y_1(j)&\ldots&\hat y_{\frac {L} 2}(j)\\
\hat y_1(j)&\hat y_2(j)&\ldots&\hat y_{\frac {L} 2+1}(j)\\
\vdots&\vdots&\vdots&\vdots\\
\hat y_{\frac {L} 2}(j)&\hat y_{\frac {L} 2+1}(j)&\ldots&\hat y_{L} (j)
\end{array}\right),
\end {equation}
has rank $m$ for $j\ne 0$ and $(m+1)/2$ for $j=0$. 
However, the matrices $\widetilde {H}(j)$ formed as in \eqref {Hj} using the noisy measurments $\tilde y_n$ will fail the rank conditions.  Cadzow's Algorithm approximates $H(j)$ via iterative changes of $\widetilde {H}(j)$  that enforce the rank and the Hankel conditions successively.

\IncMargin{1em}
\begin{algorithm}[H]\label{Cadzow denoising algorithm}
\SetKwInOut{Goal}{Goal}
\SetKwInOut{Input}{Input}
\SetKwInOut{Output}{Output}
\Goal{Denoising measurements matrix $\widetilde{Y}=\left(
\tilde y_0 \; \tilde y_1 \; \ldots \tilde y_{L}
\right)$.} 
\Input{ $\widetilde{Y}$ and $k_{\max}$ (maximal number of iterations).}
\BlankLine
{Generate the matrix $(\widetilde{Y})^\wedge$ by taking the Fourier transform on $\widetilde{Y}$.}\\
\For{$j = 0$ \KwTo $J-1$}{
{\If {$j = 0$,} {Set $r =\frac{m+1}2$}\lElse{Set $r = m$}}

{Form Hankel matrix $X = \widetilde{H}(j)$ as in \eqref{Hj} from the $j$th row of $(\widetilde{Y})^\wedge$}.
\\
\For{$k = 1$ \KwTo $k_{\max}$}{
	{Compute the SVD of $X$: $X=U\Sigma V^*$, $\Sigma = diag(\sigma_1,\ldots,\sigma_{\frac L 2 +1})$}.
	\\
	{Set $X=U diag(\sigma_1,\ldots,\sigma_r,0,\ldots,0)V^*$}.
	
	{Generate a Hankel matrix $H_{new}$ by averaging $X$ across its anti-diagonals.}
	
	{Set $X=H_{new}$.}
	
}
{Update the $j$th row of $(\widetilde{Y})^\wedge$ by the vector obtained by averaging the anti-diagonals of $X$.}
}
{Update  $\widetilde{Y}$ by taking inverse Fourier transform  
on  $(\widetilde{Y})^\wedge$.} 

\Output{Denoised data $\widetilde{Y}$.}
\caption{The pseudo-code for the Cadzow denoising method.}
\end{algorithm}
\IncMargin{1em}

For each $j\in\ZZ_J$, 
an application of the singular value decomposition (SVD) technique produces a decomposition 
$\widetilde{H}(j)=U\Sigma V^*$, where $\Sigma=diag(\sigma_1,\ldots,\sigma_{\frac{L}{2}+1})$ and $\sigma_1\geq\ldots\geq\sigma_{\frac{L}{2}+1}$. Since the rank is known to be $r_j$,  one can set    $\sigma_i=0$ for $i> r_j$ and obtain an amended matrix of singular values $\Sigma_{r_j}$. Then, one may proceed to compute  the matrix $X_{new}=U\Sigma_{r_j}V^*$ and form a new Hankel matrix   $H_{new}$  by averaging $X_{new}$ across its anti-diagonals. This procedure is  applied iteratively. After several iterations, a better approximation of the Hankel matrix $ H(j)$ is obtained and a vector of denoised data can be retrieved by applying the inverse Fourier transform.
\section {Error Analysis}\label{errant}

\subsection {Error analysis for general least squares problems}
We begin this section with the error analysis of a least squares problem that is more general than the first dynamical sampling problem.
We let $A_i\in\mathbb{C}^{m_i\times d}$, $f\in\mathbb{C}^d$, and $\tilde{y}_i=A_if+\eta_i$, where $\eta_i$ are i.i.d.~random variables with a zero mean  and  a variance matrix $\sigma^2 I$. The signal $f$ can be approximately recovered via 
\begin{eqnarray}\label{equationLSO}
f^\sharp_L=\arg\min_{g}\sum_{i=1}^{L}\|A_ig-\tilde{y}_i\|^2_2.
\end{eqnarray}
Denote the error $\epsilon_L=f^\sharp_L-f$. By \eqref{equationLSO} and the definition for $\tilde{y}_i$, it follows that
\begin{eqnarray}\label{problem5}
\epsilon_L=\arg\min_{\epsilon}\sum_{i=1}^{L}\|A_i\epsilon-\eta_i\|^2_2.
\end{eqnarray}

Let
\begin{equation}\label{equation6}
\mathcal{A}_L=\left(\begin{array}{c}
A_1\\
A_2\\
\vdots\\
A_L
\end{array}\right),
\end{equation}
and assume that for $L\geq N$, where $N$ is some fixed number,  $\mathcal{A}_L$, defined  by 
\eqref {equation6} above, has full rank. By solving problem (\ref{problem5}), we have
\begin{equation}
\epsilon_L=\left(\sum_{i=1}^{L}A_i^*A_i\right)^{-1}\sum_{i=1}^{L}A_i^*\eta_i,\text{ for all }L\geq N.
\end{equation}

The following proposition can be derived from \cite[Theorem B on p. 574]{R07}. For the convenience of the reader, however, we include the proof in the Appendix.

\begin{prop}\label{theorem2.2}
Assume that $\mathcal{A}_L$ is defined as in \eqref{equation6} and has full rank for $L\geq N$. Let $\lambda_j(L)$, $1\leq j\leq d$, denote the eigenvalues of the matrix $\mathcal{A}_L^*\mathcal{A}_L=\sum_{i=1}^{L}A_i^*A_i,1\leq j\leq d$. Then, the following holds: 
\begin{equation}
\label {estepsn}
E(\|\epsilon_L\|^2_2)=\sigma^2\sum_{j=1}^{d}1/\lambda_j(L),
\end{equation}
where $\epsilon_L$ is obtained from \eqref{problem5} and $\sigma$ is  the variance of the noise. 
\end{prop}
To study the behavior of the MSE function in \eqref{estepsn}, we recall the well-known Courant-Fischer Minimax Theorem and one of its most useful corollaries.

\begin{thm}\label{theoremCFMT}
\textbf{$\left(\text{Courant-Fischer Minimax Theorem}\right)$} Let $A$ be a $d\times d$ Hermitian matrix with eigenvalues $\lambda_1\geq\ldots\geq\lambda_k\geq\ldots\geq\lambda_d$. Then,
\[\lambda_k=\max_{U}\left\{\min_{x}\left\{\frac{x^*Ax}{x^*x}:~x\in U\text{ and }x\neq 0\right\}:~ \dim(U)=k\right\}.\]
\end{thm}
\begin{cor}\label{corollary1}
Let $A\in\mathbb{C}^{d\times d}$ and $B\in\mathbb{C}^{d\times d}$ be self-adjoint  positive semidefinite matrices. Then, $\lambda_i(A+B)\geq\lambda_i(A)$ and $\lambda_i(A+B)\geq\lambda_i(B)$.
\end{cor}
The following result is immediate from Corollary \ref{corollary1}.  \begin{prop}\label{theorem2.201}
The function $ E(\|\epsilon_L\|_2^2)$ defined by \eqref{estepsn}  is a non-negative non-increasing function of   $L$ for $L\geq N$ where $N$ is some fixed number,  such that $\mathcal{A}_N$ in 
\eqref {equation6} has full rank. Consequently, as $L$ goes to $\infty$, it converges to a non-negative constant.
\end{prop}

The goal of the following example  is to illustrate the 
above result in the context of dynamical sampling but without sub-sampling.

\begin{exmp}[Special case: no sub-sampling]\label{example2.1} 	Suppose that $A$ is a normal matrix and suppose that $A_i=A^{i-1}$ in  \eqref{problem5}. 
Because $A$ is normal, it can be written as $A=U^*DU$, where $U$ is a unitary matrix and $D$ is a diagonal matrix with the diagonal entries $s_1,s_2,\ldots,s_d$. 
Hence, $\mathcal{A}_L^*\mathcal{A}_L$ can be computed as
\begin{eqnarray}
\mathcal{A}_L^*\mathcal{A}_L&=&\sum_{k=1}^{L}(A^*)^{k-1}A^{k-1}=U^*\sum_{k=1}^{L}\left(D^*D\right)^{k-1}U.\label{equation9}
\end{eqnarray} 
Defining $\Lambda = \Lambda(L)$ by  $$\mathcal{A}_L^*\mathcal{A}_L=U^*\Lambda U,$$ and 
\[
\Lambda=\left(\begin{array}{ccc}
\lambda_1(L)\\
&\ddots&\\
&&\lambda_d(L)\\
\end{array}\right),
\]
we get from \eqref{equation9} that  $$\lambda_j(L)=\begin{cases}
\frac{1-|s_j|^{2L}}{1-|s_j|^2}, &\lvert s_j\rvert \ne 1 ;\\
L,&\lvert s_j\rvert=1.\\
\end{cases}$$ 
The error $\epsilon_L$ can be represented as
\begin{eqnarray}
\epsilon_L&=&\left(\sum_{i=1}^{L}A_i^*A_i\right)^{-1}\sum_{i=1}^{L}A_i^*\eta_i=U^*\Lambda^{-1}\sum_{i=1}^{L}(D^*)^{i-1}U\eta_i,\nonumber
\end{eqnarray}
and \eqref {estepsn} follows immediately for this special case.

To illustrate Proposition  \ref {theorem2.201},  note that, when $|s_j|<1$, the expression $\frac 1 {\lambda_j(L)}=\frac{1-|s_j|^2}{1-|s_j|^{2L}}$  decreases  and  converges to  $1-|s_j|^2$ as $n$ increases and tends to $ \infty$. 

When $|s_j|=1$, then $\frac 1 {\lambda_j(L)}=\frac{1}{L}$ which  decreases as $L$  increases.

When $|s_j|>1$, $\frac 1 {\lambda_j(L)}=\frac{1-|s_j|^2}{1-|s_j|^{2L}}$   decreases (as $L$  increases) and  converges to  $0$ as $L \to \infty$. 

Thus, in all  three cases   the function $E(\|\epsilon_L\|_2^2)$ is decreasing as $L$ increases. In addition,
\[E(\|\epsilon_L\|_2^2)\rightarrow \sigma^2\sum_{\substack{1\leq j\leq d\\|s_j|<1}}(1-|s_j|^2), \text{ as } L\rightarrow\infty.\]
\end{exmp}
\subsection { Error analysis for dynamical sampling}
To derive a similar result for dynamical sampling, we replace the general operator $A_i$ in \eqref{problem5} with the $i-1$ power $A^{i-1}$ of a matrix $A$ followed by a subsampling matrix $S_{\Omega} $, i.e., we let $A_i=S_{\Omega} (A^{i-1})$.
By Propositions \ref{theorem2.2} and \ref{theorem2.201}, and using the fact that  $S_{\Omega}^*S_{\Omega} = S_{\Omega} $,   the following assertions hold.

\begin {thm}
\label {AympConst}
Let $\lambda_j(L)$ denote the $j$-th eigenvalue of the matrix $\sum_{i=0}^{L-1}(A^*)^i S_{\Omega}A^i$. Then
\[E(\|\epsilon_L\|_2^2)=\sigma^2\sum_{j=1}^{d}1/\lambda_j(L),\] 
is non-increasing as a function of $L$. Hence, it converges to some constant as $L \to \infty$.
\end {thm}

\begin {rem} \label {ChoiceOmega}
The theorem above shows how the mean squared error depends on  $\Omega$, $A$ and $L$. However, for a given $A$, not all choices of $\Omega$ are allowable: there are necessary and sufficient conditions on the choice of $\Omega$  that will allow us to reconstruct $f$ by solving \eqref {equationLSO}  when $A_i=S_{\Omega} (A^{i-1})$ and no noise is present \cite{ACMT17} (i.e., $A_1$ is full rank and $\lambda_j(L)>0$ for all $j,L$) . 
\end{rem}

\section {Numerical Results}\label{test}
\subsection {Error Analysis}\label{sim1}
In this section, we  illustrate the performance of the least squares based method for signal recovery (i.e., Algorithm \ref{algorithm iterative}) in the case when the dynamical samples are corrupted by noise. We describe the numerical simulations that we conducted using synthetic data and examine the behavior of  $E(\|\epsilon_L\|_2^2/\sigma^2)$  as a function of the number of time levels $L$. 
\begin{figure}[tbph]
\begin{subfigure}[b]{0.49\textwidth}
\includegraphics[width=\textwidth]{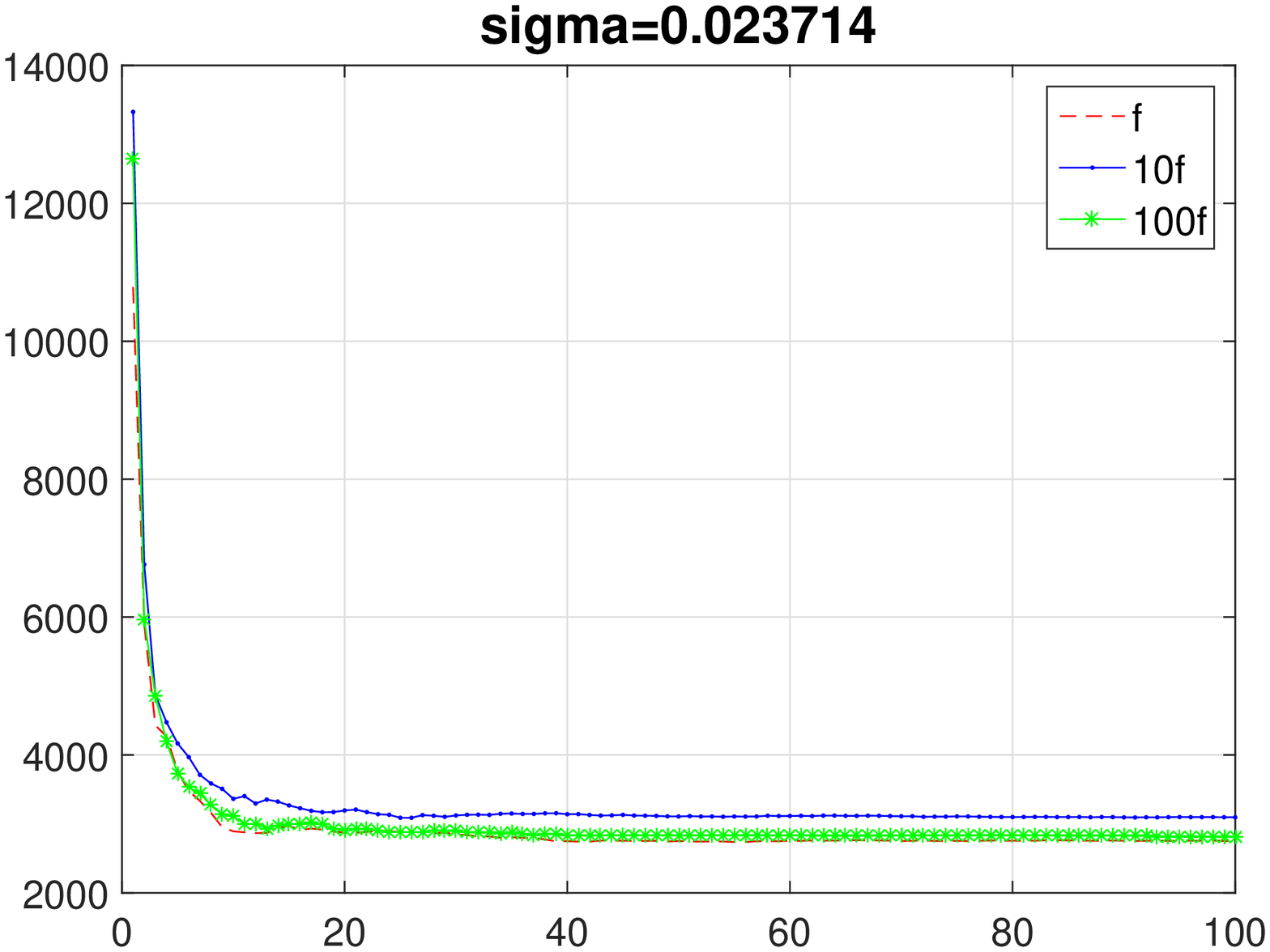}
\caption{}
\label{fig:1.1}
\end{subfigure}
\begin{subfigure}[b]{0.49\textwidth}
\includegraphics[width=\textwidth]{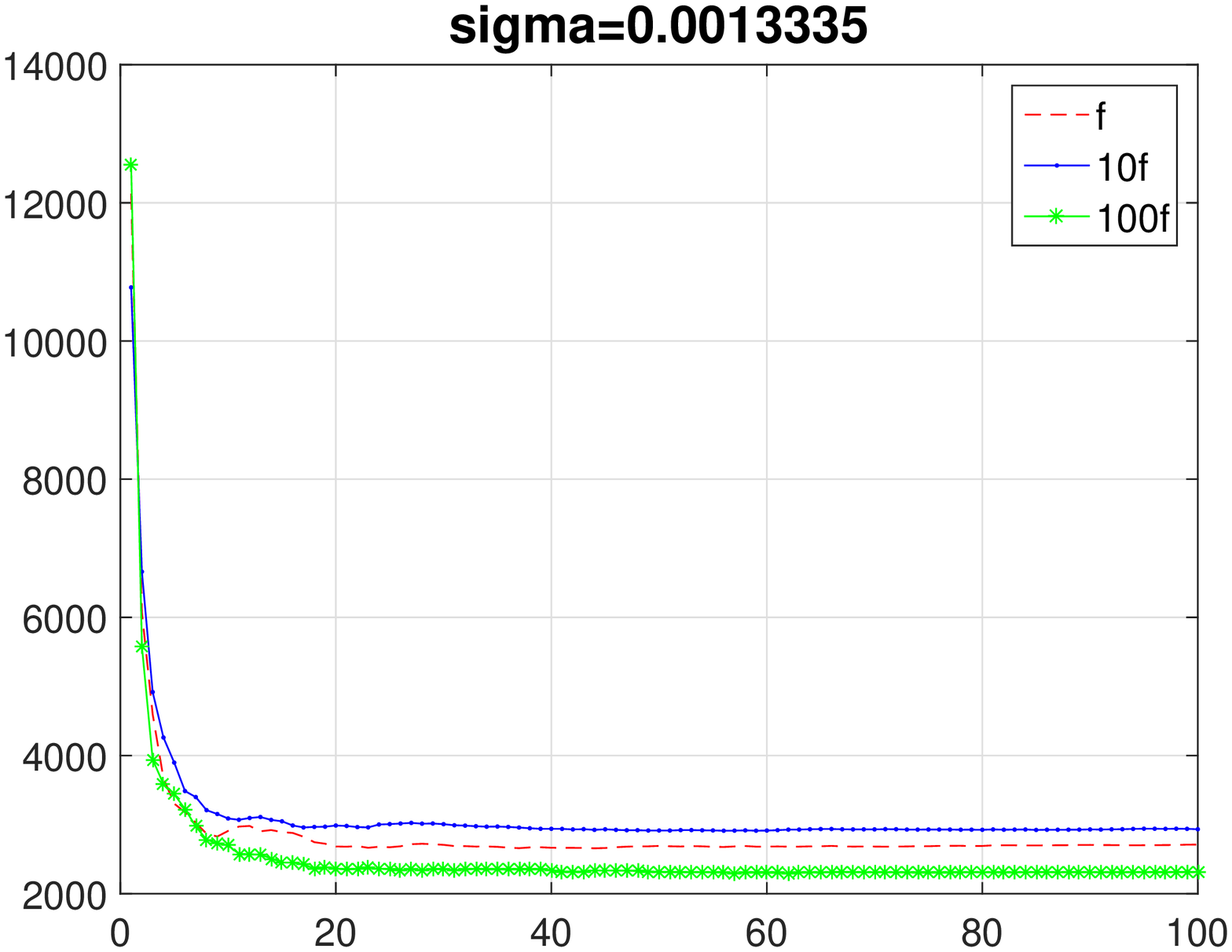}
\caption{}
\label{fig:1.2}
\end{subfigure}
\caption{The behavior of $E(\|\epsilon_L\|_2^2/\sigma^2)$ for different signals. Here, the curves in (\ref{fig:1.1}) and (\ref{fig:1.2}) show the results with noise variances $2.3714\times 10^{-2}$ and $1.3335\times 10^{-3}$, respectively. The signal $f$ is randomly generated with norm 2.2914. Three signals $f$, $10f$, and $100f$ are used for the simulations, where  $x$-axis stands for  the time levels and $y$-axis represents the value of  $E(\|\epsilon_L\|_2^2/\sigma^2)$.}
\label{fig:1}
\end{figure}

\begin{figure}
\centering
\begin{subfigure}[b]{0.49\textwidth}
\includegraphics[width=\textwidth]{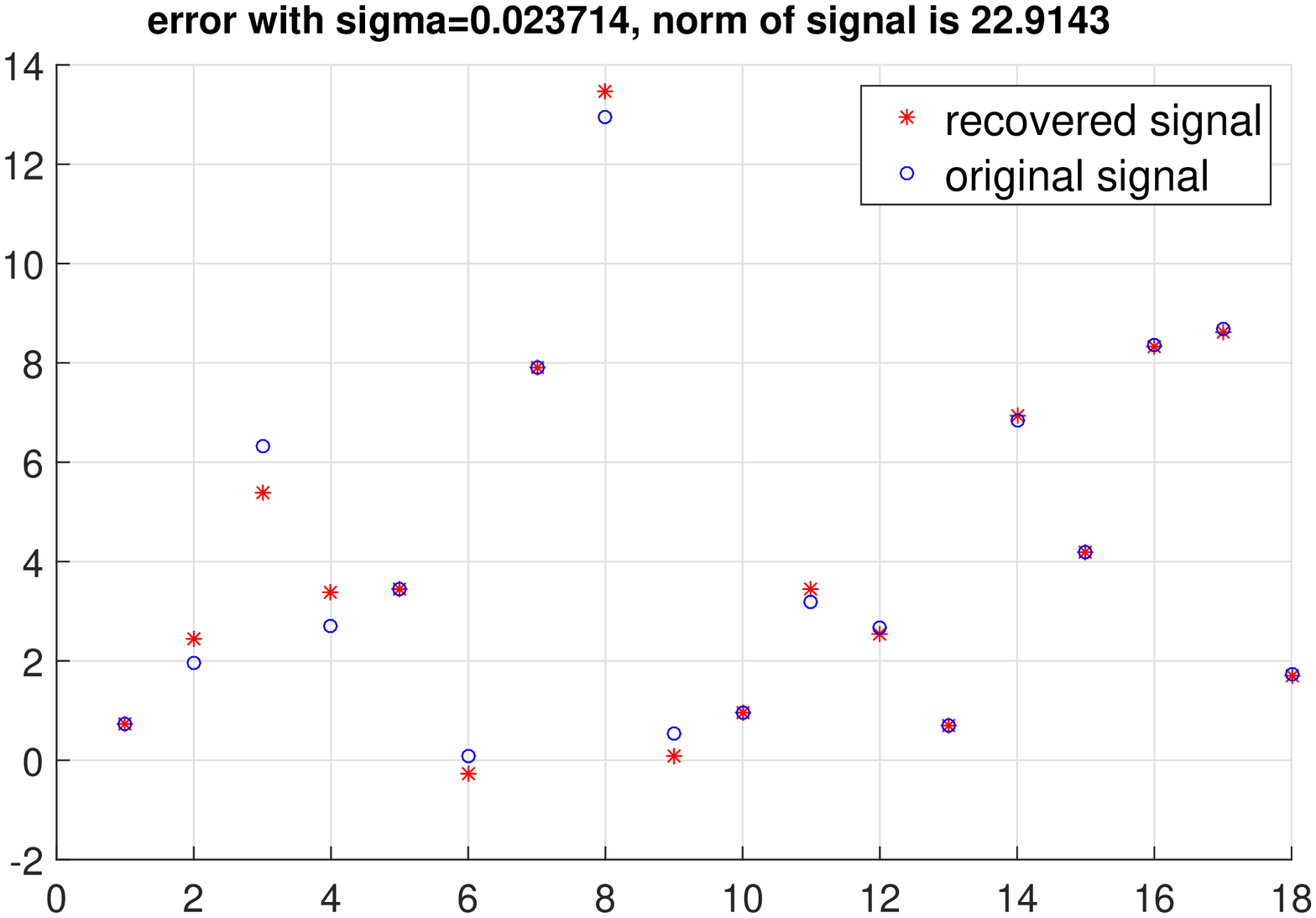}
\caption{}		\label{fig:2.1}
\end{subfigure}
\begin{subfigure}[b]{0.49\textwidth}
\includegraphics[width=\textwidth]{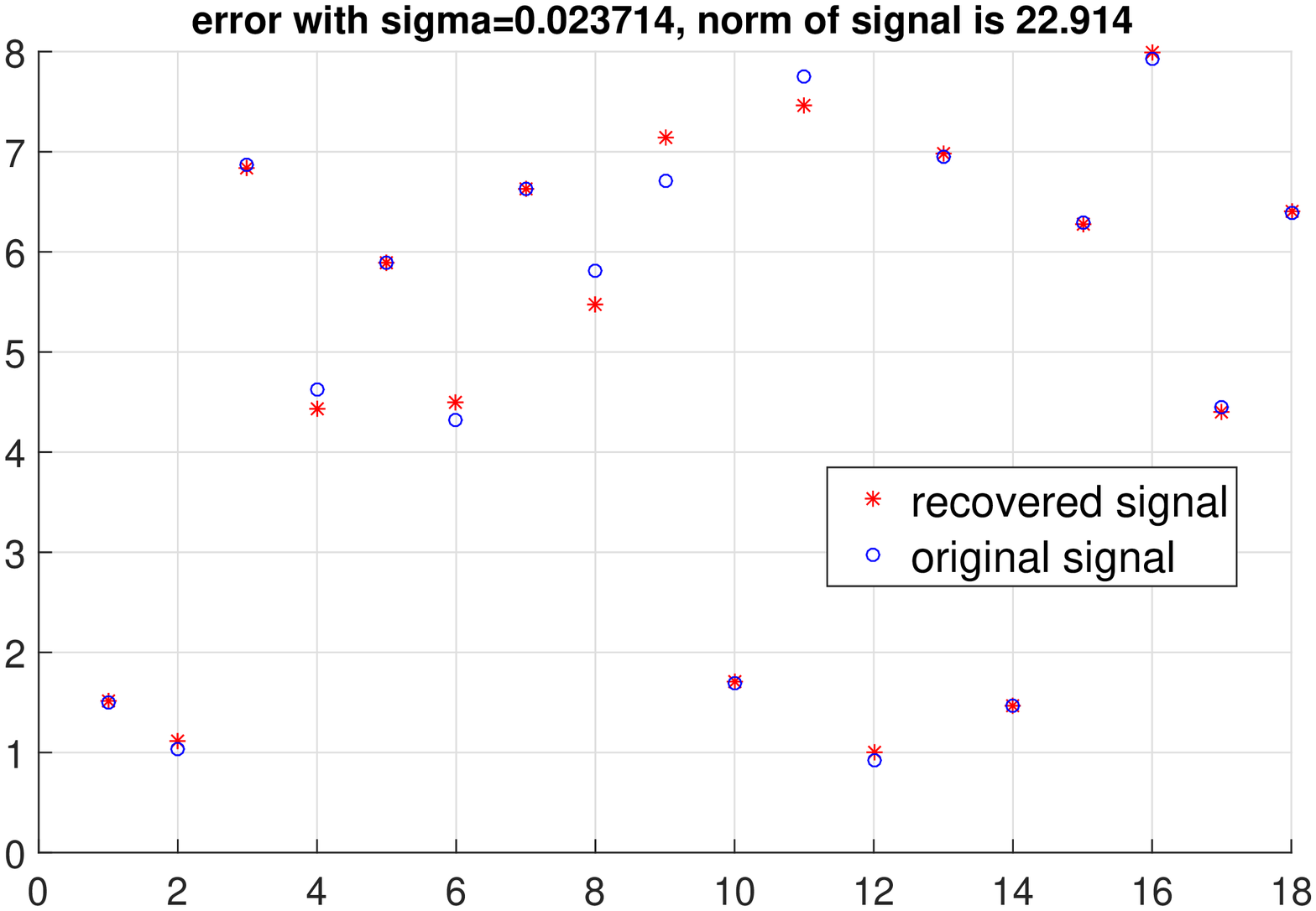}
\caption{}		\label{fig:2.2}
\end{subfigure}
\\
\begin{subfigure}[b]{0.49\textwidth}
\includegraphics[width=\textwidth]{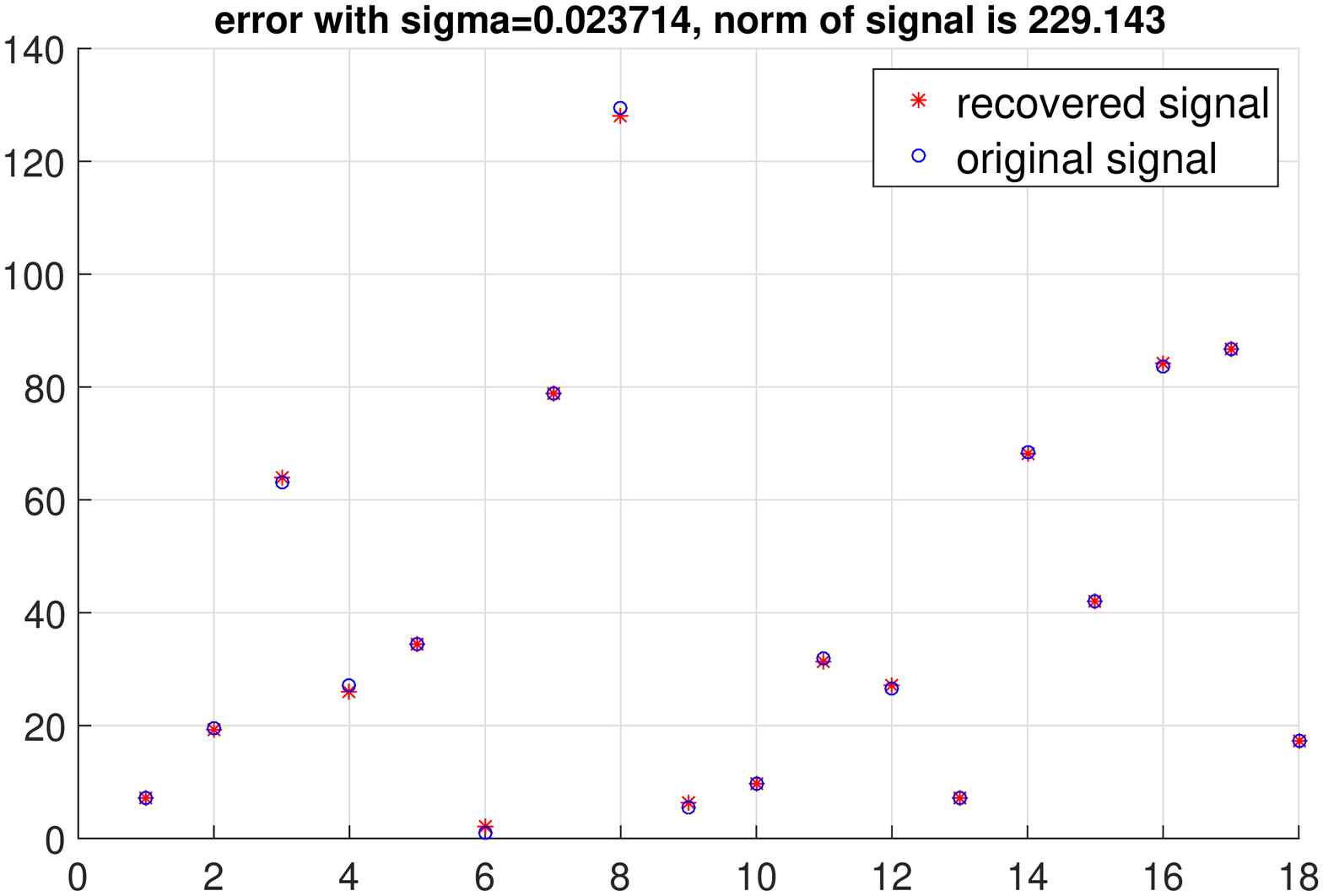}
\caption{}		\label{fig:2.3}
\end{subfigure}
\begin{subfigure}[b]{0.49\textwidth}
\includegraphics[width=\textwidth]{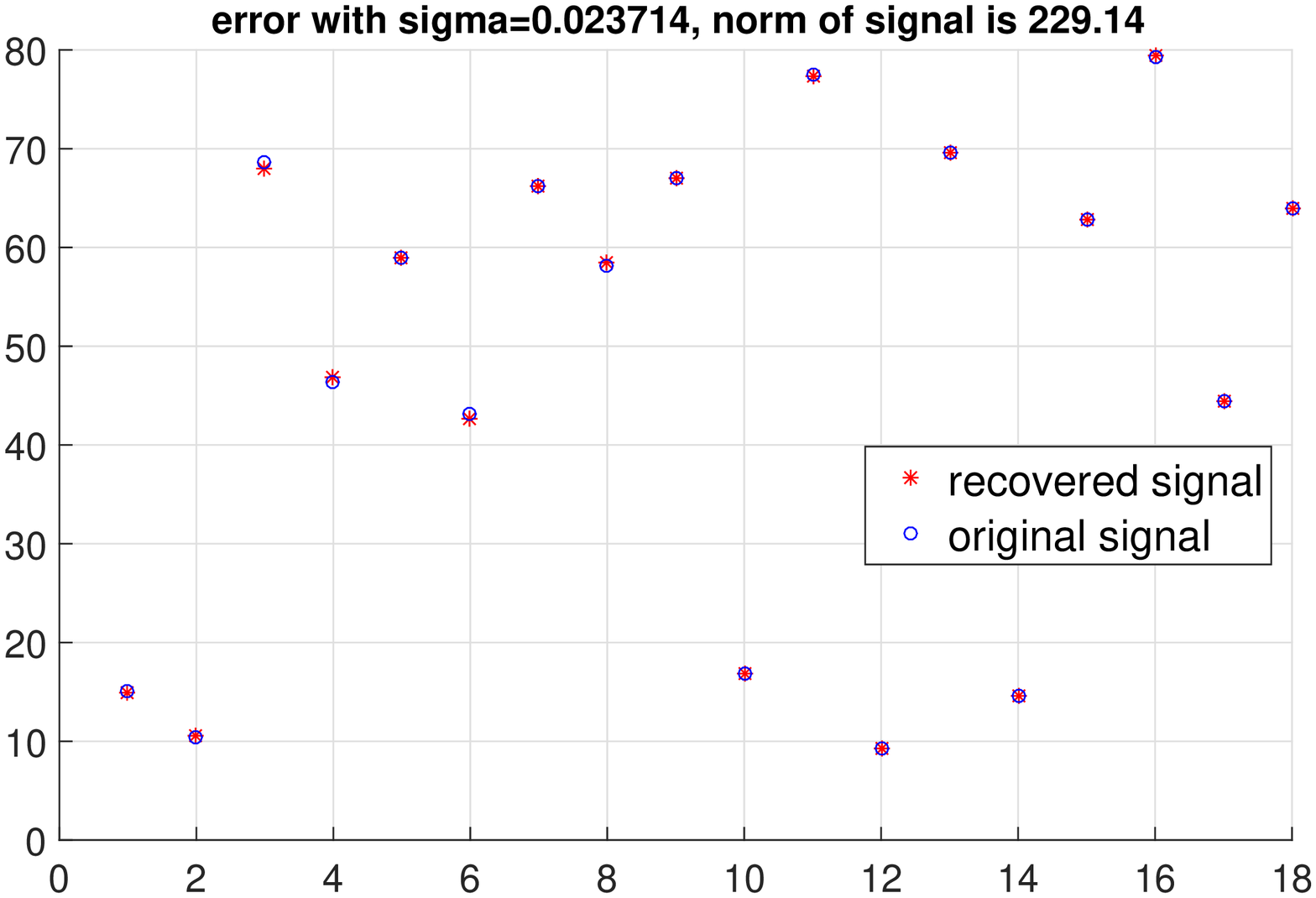}
\caption{}		\label{fig:2.4}
\end{subfigure}

\caption{The original signals and the reconstructed signals are represented by blue circles and red stars, respectively. The  signals in \ref{fig:2.3} and \ref{fig:2.4} are obtained from the signals in \ref{fig:2.1} and \ref{fig:2.2}, respectively, by multiplying  by 10. The norm of the original signal in \ref{fig:2.1} equals the norm of the original signal in \ref{fig:2.2}, the same is true for the signals in \ref{fig:2.3} and \ref{fig:2.4}.}
\label{fig:2}
\end{figure}

To obtain synthetic data for the simulation, we use   a  random  signal $f\in\ell^2(\mathbb{Z}_{18})$ and a convolution operator  $Af=a\ast f$, determined by a real symmetric vector $a$ with non-zero components given by $(\frac{1}{8}, \frac{1}{2},1  ,\frac{1}{2} , \frac{1}{8})$, i.e., $A\in\mathbb{R}^{18\times 18}$ is a circulant matrix with the first row  $\left(1,1/2,1/8,0,\ldots,0,1/8,1/2\right).$
We generate the signals $f_i=A^if$ at time levels $i=0,1,\ldots, L$.  The non-uniform locations $\Omega=\{1,5,7,10,13,15,18\}$ are chosen to generate the samples $\{f_i(j): j \in \Omega\}$. Independent and identically distributed Gaussian  noise with zero mean is then  added to the samples to obtain a set of  noisy data $\{f_i(j)+\eta_i(j): j \in \Omega\}$.

Figure \ref{fig:1} shows the relationship between $E(\|\epsilon_L\|_2^2/\sigma^2)$ and the time levels, where $\epsilon_L$ is defined by \eqref{problem5}. For each $L$, the simulation was repeated  100 times with the same distribution of noise, and $E(\|\epsilon_L\|_2^2/\sigma^2)$ was estimated by averaging the 100 values of  $\|\epsilon_L\|_2^2/\sigma^2$. 
Figure \ref{fig:1.1} shows how $E(\|\epsilon_L\|_2^2/\sigma^2)$ changes as $L$ varies for three different signals: $f$, $10f$, and $100f$, where  the noise variance is $\sigma=2.3714 \times 10^{-2}$ and  the 2-norm of $f$ approximately equals $2.2914$. The graph of $10f$ is given in Figure \ref{fig:2.1}.  Figure \ref{fig:1.2} shows the behavior of $E(\|\epsilon_L\|_2^2/\sigma^2)$ for the same signals  as in Figure \ref{fig:2.1}, where  the noise variance is $\sigma=1.3335\times10^{-3}$.
As shown in Figure \ref{fig:1}, $E(\|\epsilon_L\|_2^2/\sigma^2)$ decreases  as $n$ increases and approaches the constant  predicted by Theorem \ref  {AympConst}. 
\begin{figure}[tbph]
\begin{subfigure}[b]{0.49\textwidth}
\includegraphics[width=\textwidth]{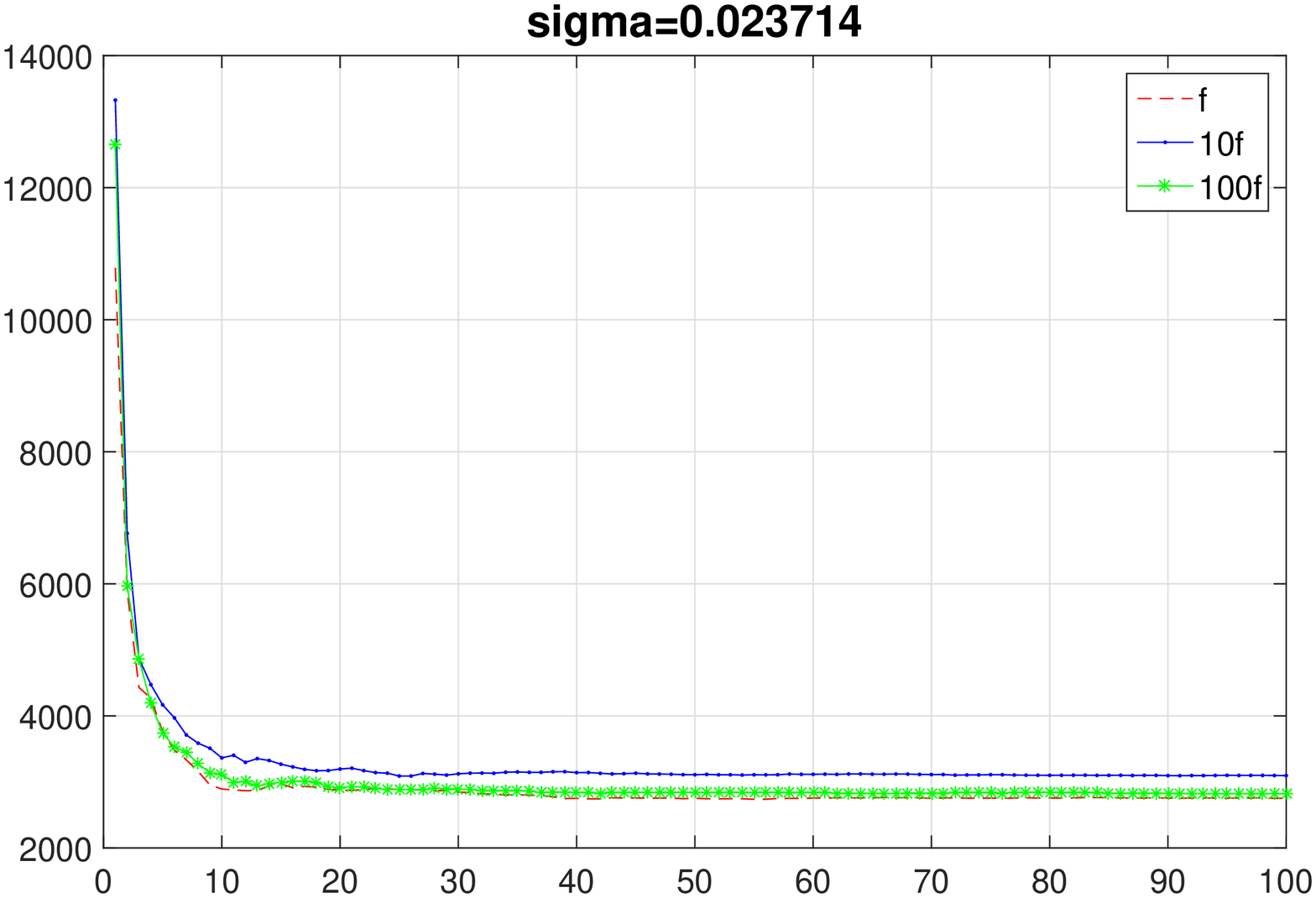}
\caption{}
\label{fig:3.1}
\end{subfigure}
\begin{subfigure}[b]{0.49\textwidth}
\includegraphics[width=\textwidth]{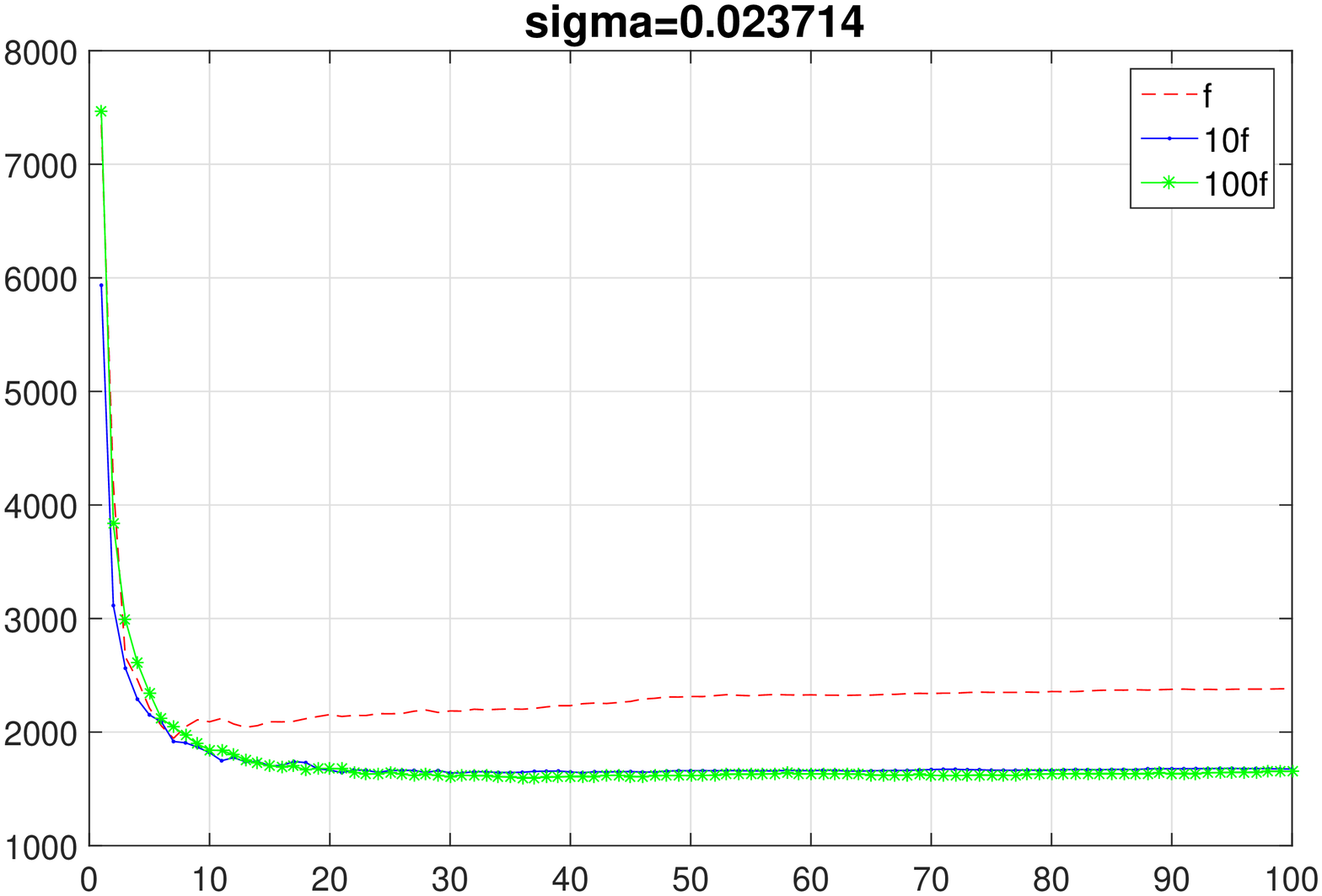}
\caption{}
\label{fig:3.2}
\end{subfigure}
\caption{ The behavior of  $E(\|\epsilon_L\|_2^2/\sigma^2)$ are shown in \eqref{fig:3.1} and \eqref{fig:3.2} for the sparsely supported signals without and with applying the threshold method, respectively, where the samples are corrupted by  Gaussian noise with zero mean  and standard deviation $2.3714\times 10^{-2}$.}
\label{fig:3}
\end{figure}

Figure \ref{fig:2}  depicts the graphs  of the reconstructed signals  and the original signals $10f$, $100f$, $10g$, and $100g$ in Figures \ref{fig:2.1}, \ref{fig:2.3}, \ref{fig:2.2}, and \ref{fig:2.4}, respectively, where $f, g\in\mathbb{R}^{18}$ are randomly generated and scaled to the norm $2.2914$. The reconstructed signals from the noisy data and the original signals are shown  in Figure \ref{fig:2.1} for $10f$,  in Figure \ref{fig:2.2} for $10g$, in Figure \ref{fig:2.3} for $100f$, and in Figure \ref{fig:2.4} for $100g$, respectively. The noisy data are corrupted  by  Gaussian noise with zero mean  and standard deviation $2.3714\times 10^{-2}$.   As displayed in Figure \ref{fig:2}, while the reconstructed signals in Figures \ref{fig:2.1} and \ref{fig:2.2} are clearly different from the original signals, it is hard to distinguish the reconstructed signals from the original signals in Figures \ref{fig:2.3} and \ref{fig:2.4} because the reconstructed signals are very close to the  original signals.

\begin{figure}[tbph]
\begin{subfigure}[b]{0.49\textwidth}
\includegraphics[width=\textwidth]{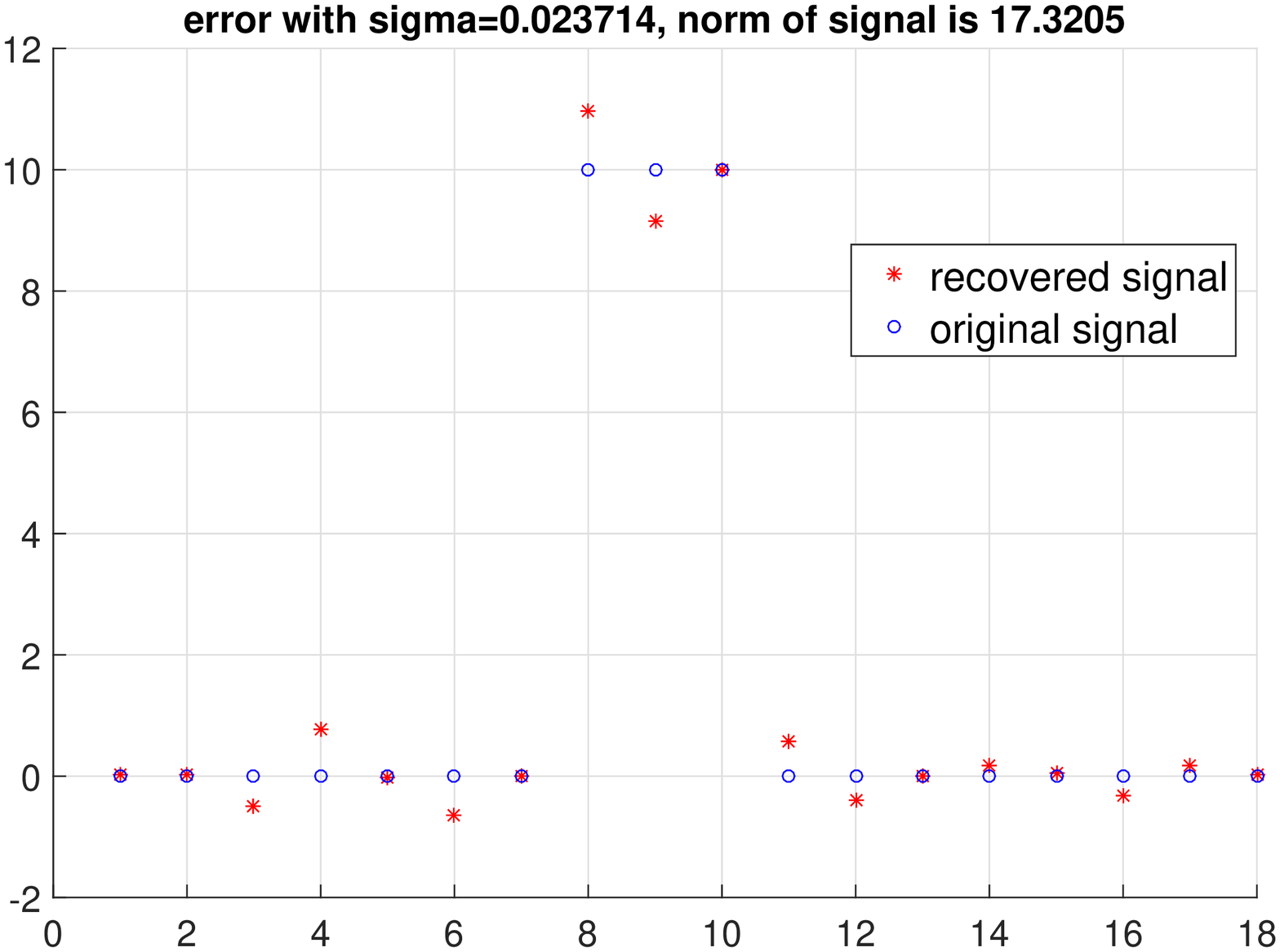}
\caption{}
\label{fig:4.1}
\end{subfigure}
\begin{subfigure}[b]{0.49\textwidth}
\includegraphics[width=\textwidth]{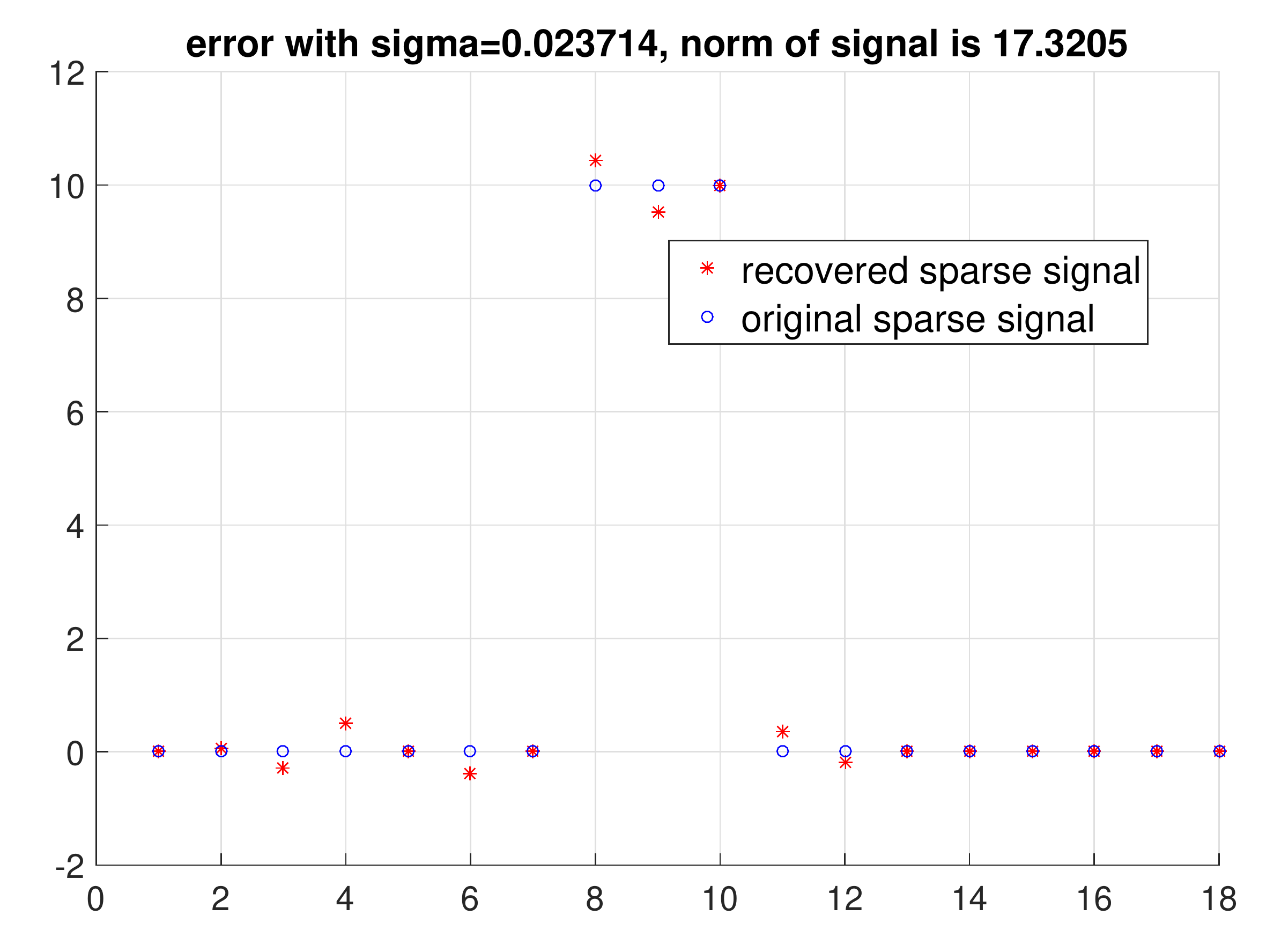}
\caption{}
\label{fig:4.2}
\end{subfigure}
\caption{A comparison of the reconstruction results before and after applying the  threshold method. (\ref{fig:4.1}) and (\ref{fig:4.2}) show the reconstruction results before and after applying the threshold method, respectively. In  (\ref{fig:4.1}) and (\ref{fig:4.2}), the original signals have the same sparse support $\{8,9,10\}$. The samples are corrupted by the independent Gaussian noise with mean $0$ and standard deviation $2.3714\times 10^{-2}$. }
\label{fig:4}
\end{figure}


Figures \ref{fig:3} and \ref{fig:4} are  simulation results for the sparsely supported signals. For these special signals, a threshold method  \cite{R07} is introduced for the samples and reconstructed signals. The  method is implemented as follows. Let the threshold $T$ be $2\sigma$, 
let  $\tilde y$ denote the sample vector, and let $f_L^\sharp$ be  the reconstructed signal. If $\lvert \tilde y(i)\rvert\le T$, we set  $\tilde y(i)=0$, where $\tilde y(i)$  is the $i$-th component of $\tilde y$. Similarly, if $\lvert f_L^\sharp(i)\rvert \le T$, we set   $f_L^\sharp(i)=0$.    Then the reconstruction results before and after applying 
the threshold method  are compared. 
Figures \ref{fig:3} and \ref{fig:4}  illustrate the behavior of $E(\|\epsilon_L\|_2^2/\sigma^2)$ and the reconstructed signals before and after applying the threshold method, respectively.
In the simulation, the samples are corrupted by  Gaussian noise with zero mean  and standard deviation $2.3714\times 10^{-2}$. A sparsely supported signal $f\in \mathbb{R}^{18}$  is generated with support in the locations $\{8,9,10\}$ with $f(8)=f(9)=f(10)= 1$. The MSE $E(\|\epsilon_L\|_2^2/\sigma^2)$ are estimated   for signals $f$, $10f$, and $100f$ separately. As shown in Figure \ref{fig:3}, for $n$ sufficiently large,  $E(\|\epsilon_L\|_2^2/\sigma^2)$ is about $20\%$ smaller after the threshold method is applied to the samples and reconstructed signals. Figure \ref{fig:4} shows the graphs of the original signal $10f$ and the reconstructed signal, which suggests that  the signal reconstructed by applying the threshold method is more accurate than the one reconstructed without applying the threshold method in the locations outside the support of the original signal. These observations suggest that  the threshold method can reduce $E(\|\epsilon_L\|_2^2/\sigma^2)$ by improving  the accuracy of the zero sets.

\subsection{Cadzow Denoising} \label{cd}
In this section, we describe the impact of the Cadzow denoising technique described in Section \ref {noise_red} on dynamical sampling  using synthetic data. 

\subsubsection{Denoising of the sampled data}\label{section sample reconstruction} We use a symmetric convolution operator 
$A$  with eigenvalues $\{1/8,1/4,3/8,1/2,5/8,3/4,7/8,1\}$. We let   $A$ act on the normalized randomly generated signal $f=(0.2931,$ 0.3258, 0.04568, 0.3286, 0.2275,  0.0351, 0.1002, 0.1967, 0.3444, 0.34710, 0.0567, 0.3492, 0.3443, 0.1746, $0.2879)^{T}$ iteratively for 100 times. The iterated signals are stored in a matrix $\Pi$ as
$$ \Pi=\left(
f\; Af \;  A^2f\; \ldots\; A^{100}f
\right)=(f_0 \; f_1\; f_2 \ldots \; f_{100})$$ where $A^{k}f$ is a column vector for each $0\leq k\leq 100$ (see   \eqref {equation1}). At each time level, the generated signals are perturbed by i.i.d.~Gaussian noise with zero mean   and standard deviation $\sigma\in\{10^{-2},10^{-3},10^{-4},10^{-5}\}$; the noisy signals are denoted by $$\widetilde{\Pi}=\Pi+H,$$ where $H_{i,j}\sim N(0,\sigma^2)$ and every two entries of $H$ are independent (see \eqref{NoisyDySamp}).  

The samples are taken uniformly on $3\ZZ_{15}$ (i.e., $m=3$), specifically  at locations $\Omega=\{1,4,7,10,13\}$. The Cadzow algorithm (Algorithm \ref {Cadzow denoising algorithm}) is applied to the data  $\widetilde Y=S_m\tilde{\Pi}$ where $S_m$  is defined in the first paragraph of Section \ref {frec}.   The denoised data are denoted by   $Z$ which is compared to $S_m\Pi$ directly by computing  \begin{equation}\label{dnoeq}
\frac{\|Z-S_{m}\Pi\|}{\|S_m\Pi\|}.
\end{equation}
In addition,   the relative difference between the noisy data $S_m\widetilde{\Pi}$ and $S_m\Pi$ is  computed as
\begin{equation}\label{noeq}
\frac{\|S_m\widetilde{\Pi}-S_m\Pi\|}{\|S_m\Pi\|}. \end{equation}

The same process is repeated for 80 times (with the same $\Pi$ and different $H$). The numerical results are obtained by  averaging   the 80 values of  
\eqref{dnoeq}  and \eqref{noeq}, respectively.

\begin{figure}[tbph]
\begin{center}
\includegraphics[width=0.6\textwidth]{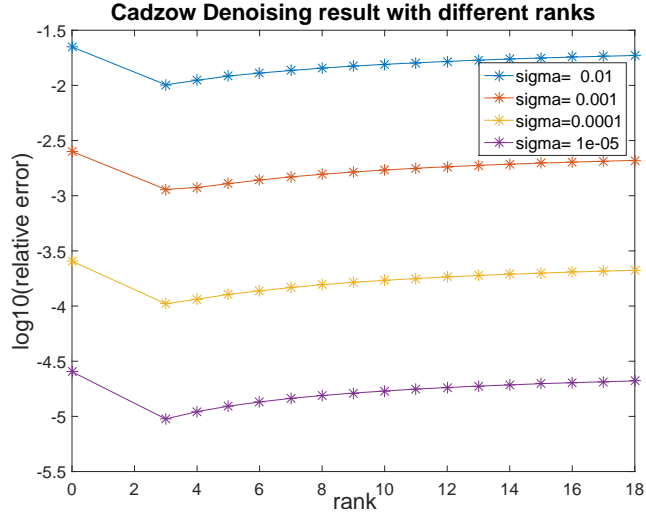}
\end{center}
\caption{The relative errors using  Cadzow denoising method. The vertical axis  represents $\log_{10}$   of averaged \eqref{dnoeq} when the threshold rank  in the Cadzow denoising Algorithm \ref {Cadzow denoising algorithm} is greater than or equal to $ 3$. When the threshold rank is $0$  \eqref{noeq} is used instead of \eqref{dnoeq}. }\label{figurecadtest}
\end{figure}

The  simulation results are  shown in Figure \ref{figurecadtest}. For $rank\ge 3$, the horizontal values depict  the threshold ranks in the Cadzow algorithm.  The corresponding vertical values are   $\log_{10}$ of the  values of  \eqref {dnoeq} averaged over 80 repetitions. When $rank=0$, \eqref{noeq} is used instead of \eqref {dnoeq}.   As  shown in  Figure \ref{figurecadtest}, the Cadzow denoising technique works best for noise  reduction when the rank of the Hankel matrix is chosen to be $3$, which is consistent with the theory  described  in Section \ref{Cadzow denoising theory}.

\begin{figure}
\begin{center}
\includegraphics[width=0.7\textwidth]{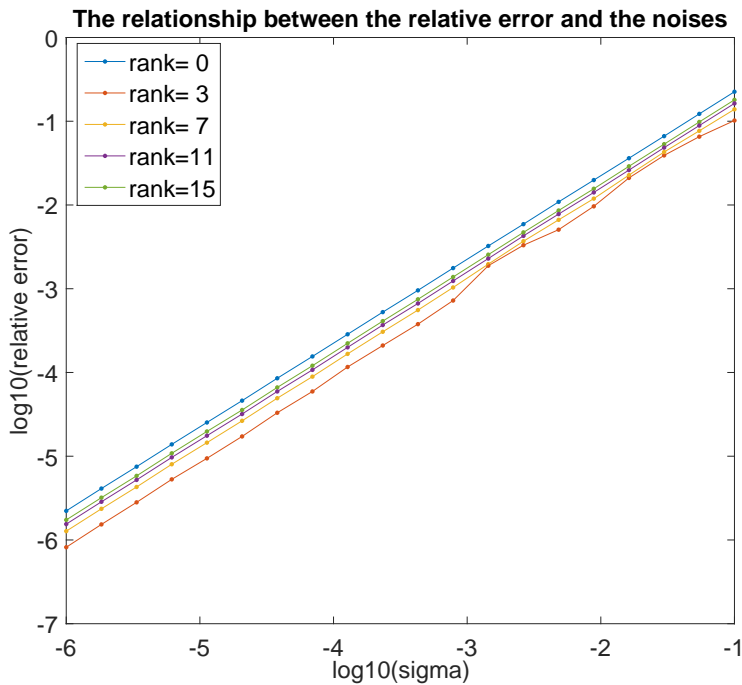}
\end{center}
\caption{The relation between the relative errors and the noise standard deviations using the  Cadzow denoising method with different threshold ranks. The  curves labeled  $``\text{rank}=r"$, $r=3,7\dots$ reflect the relationship between  $\log_{10}$  of the averaged relative errors  and $\log_{10}$ of the noise standard deviations, where the relative errors are represented in \eqref{dnoeq} for $rank\geq 3$ and in \eqref{noeq} for $rank=0$.} \label{figureerrornoise}
\end{figure}

In Figure \ref{figureerrornoise}, the curve labeled $``\text{rank}=0"$ shows the relationship between $\log_{10}$ of the averaged \eqref{noeq}  and  $\log_{10}$ of the noise standard deviations, while the  curves labeled  as $``\text{rank}=r"$ for $r=3,7,11,15$ show the relationship between  $\log_{10}$  of the averaged \eqref{dnoeq} and $\log_{10}$ of the noise standard deviation. As displayed in Figure \ref{figureerrornoise}, the curves are almost linear.  For fixed noise standard deviation, the figure shows that, as predicted by the theory described in Section \ref {Cadzow denoising theory},  the best denoising happens when $``\text{rank} =3"$ since the sub-sampling is $3$. 

\subsubsection{Spectrum Reconstruction of the Convolution Operator}
\begin{figure}[tbph]
\centering
\begin{subfigure}[b]{0.49\textwidth}
\includegraphics[width=\textwidth]{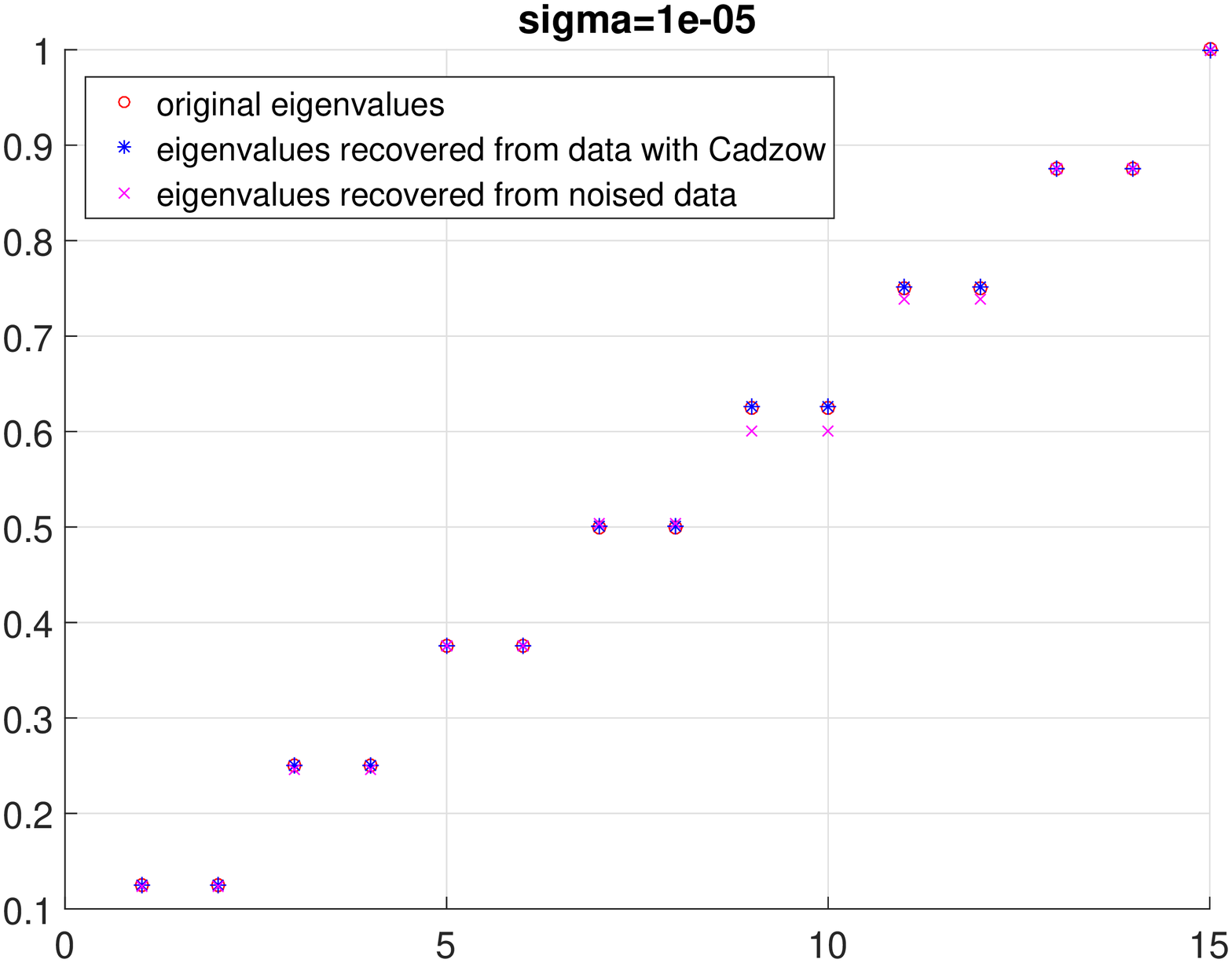}
\caption{}
\label{fig:4.2.1.1a}
\end{subfigure}
\begin{subfigure}[b]{0.49\textwidth}
\includegraphics[width=\textwidth]{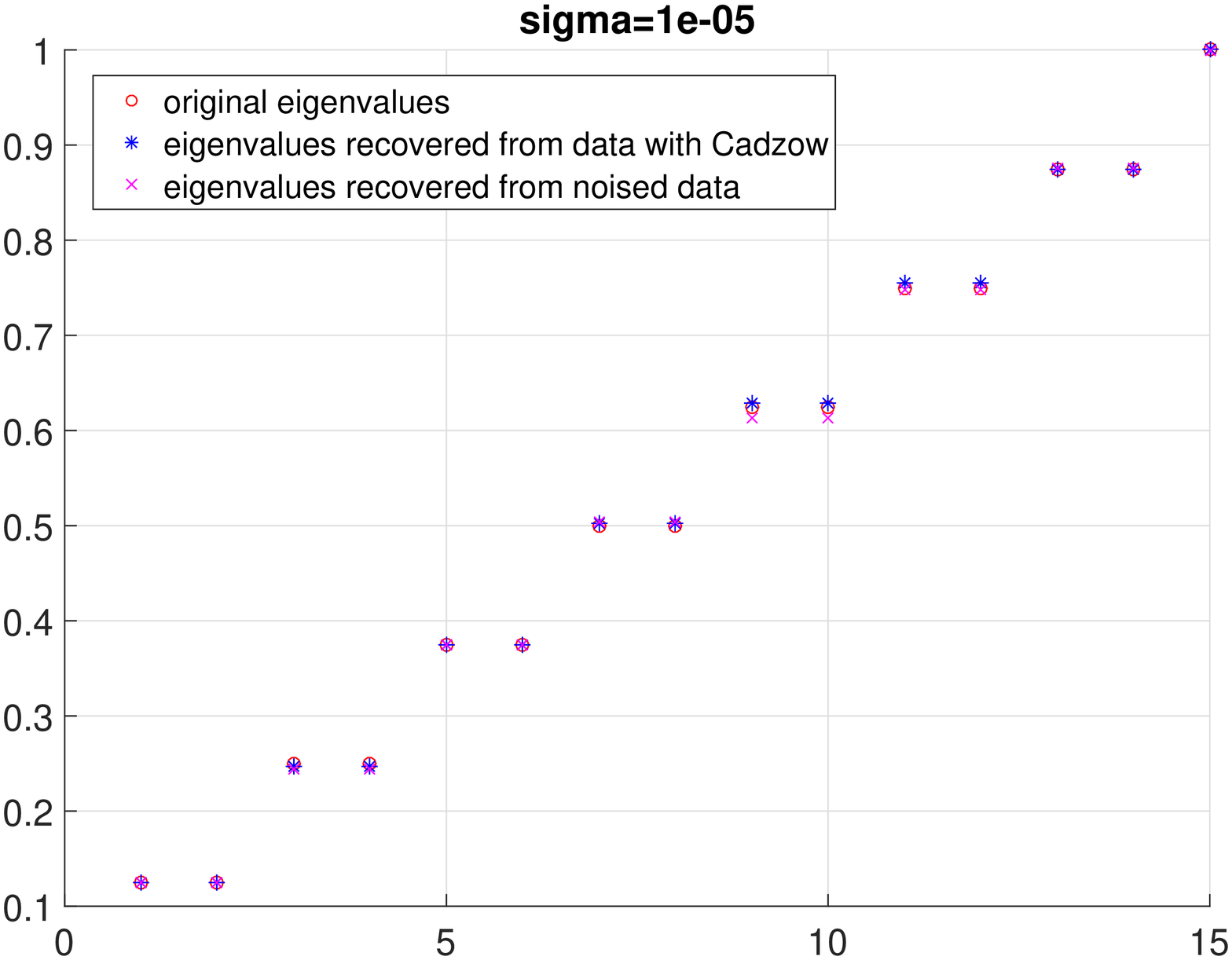}
\caption{}
\label{fig:4.2.1.1b}
\end{subfigure}
\begin{subfigure}[b]{0.49\textwidth}
\includegraphics[width=\textwidth]{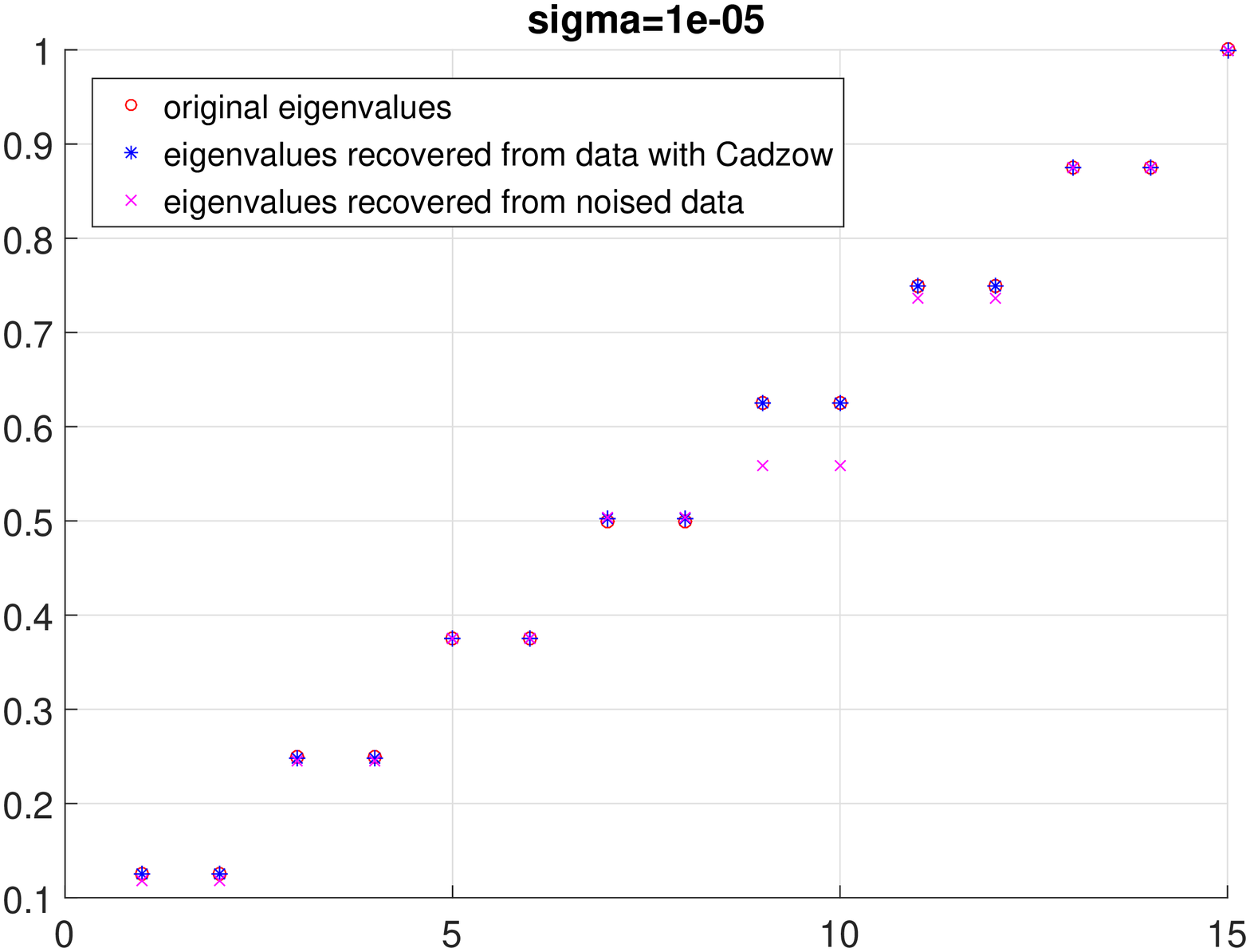}
\caption{}
\label{fig:4.2.1.1c}
\end{subfigure}
\caption{A comparison of  the spectrum reconstruction with 
and without the Cadzow denoising technique  for $\sigma=10^{-5}$.}
\label{fig:4.2.1.1}
\end{figure}

\begin{figure}[tbph]
\centering
\begin{subfigure}[b]{0.49\textwidth}
\includegraphics[width=\textwidth]{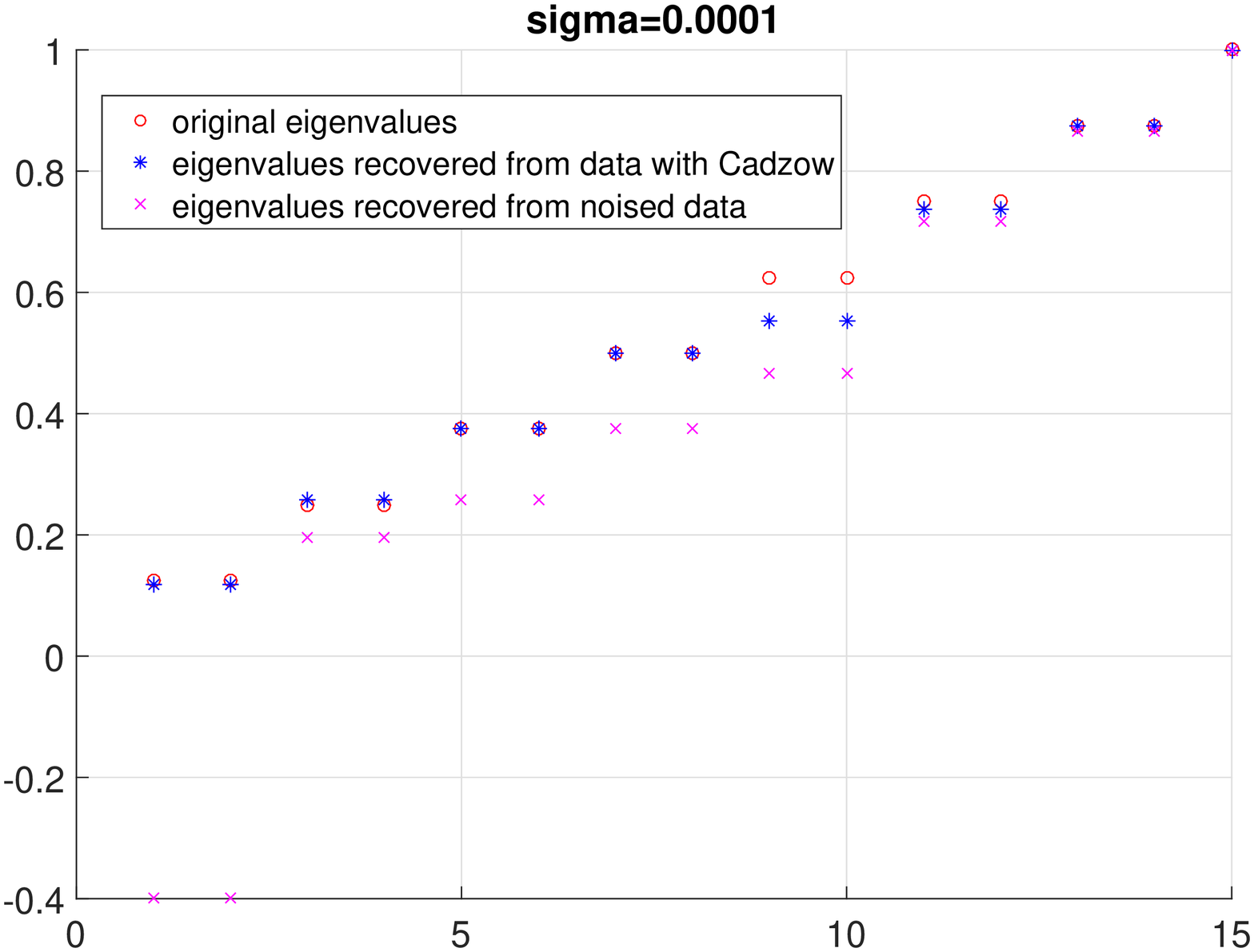}
\caption{}
\end{subfigure}
\begin{subfigure}[b]{0.49\textwidth}
\includegraphics[width=\textwidth]{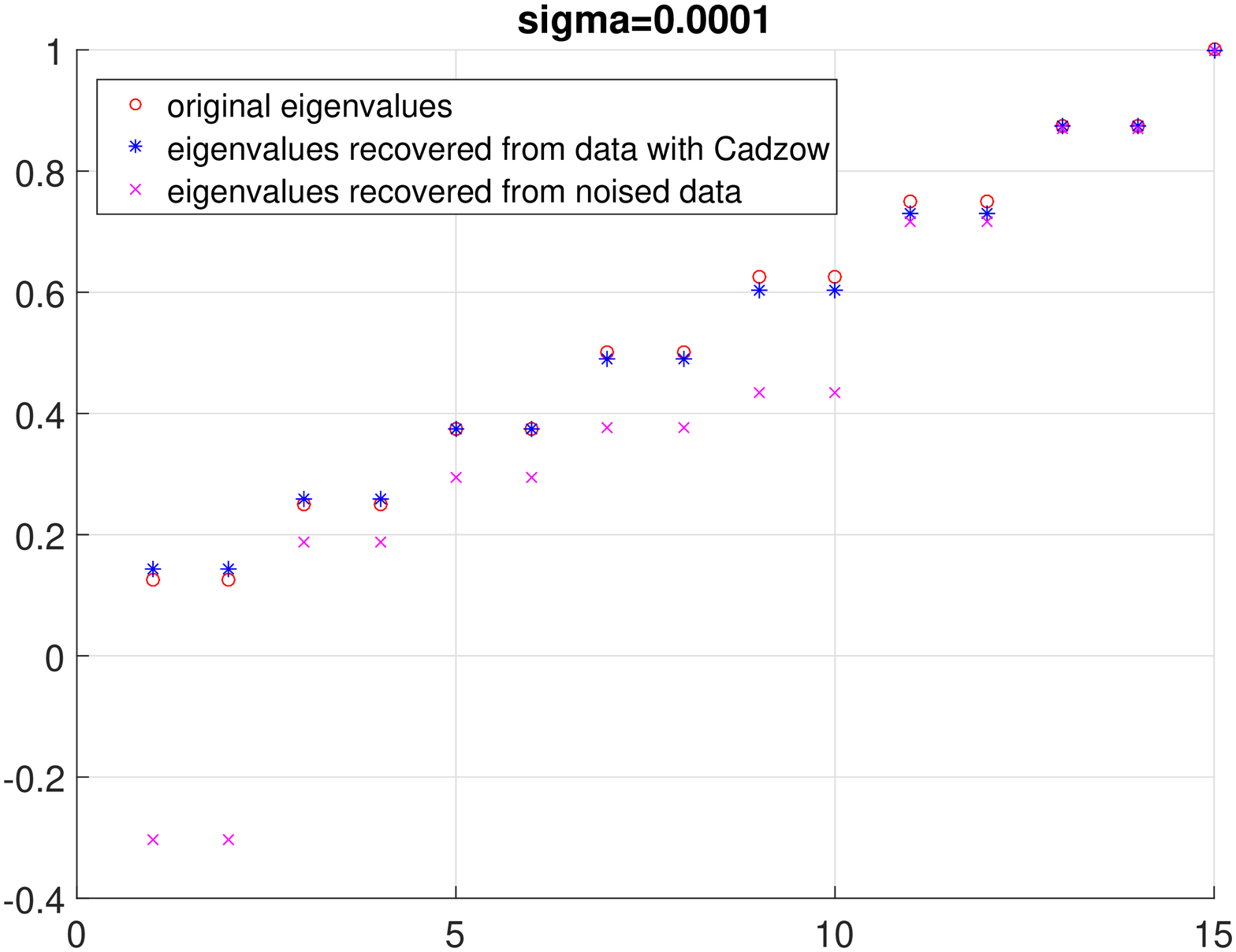}
\caption{}
\end{subfigure}
\begin{subfigure}[b]{0.49\textwidth}
\includegraphics[width=\textwidth]{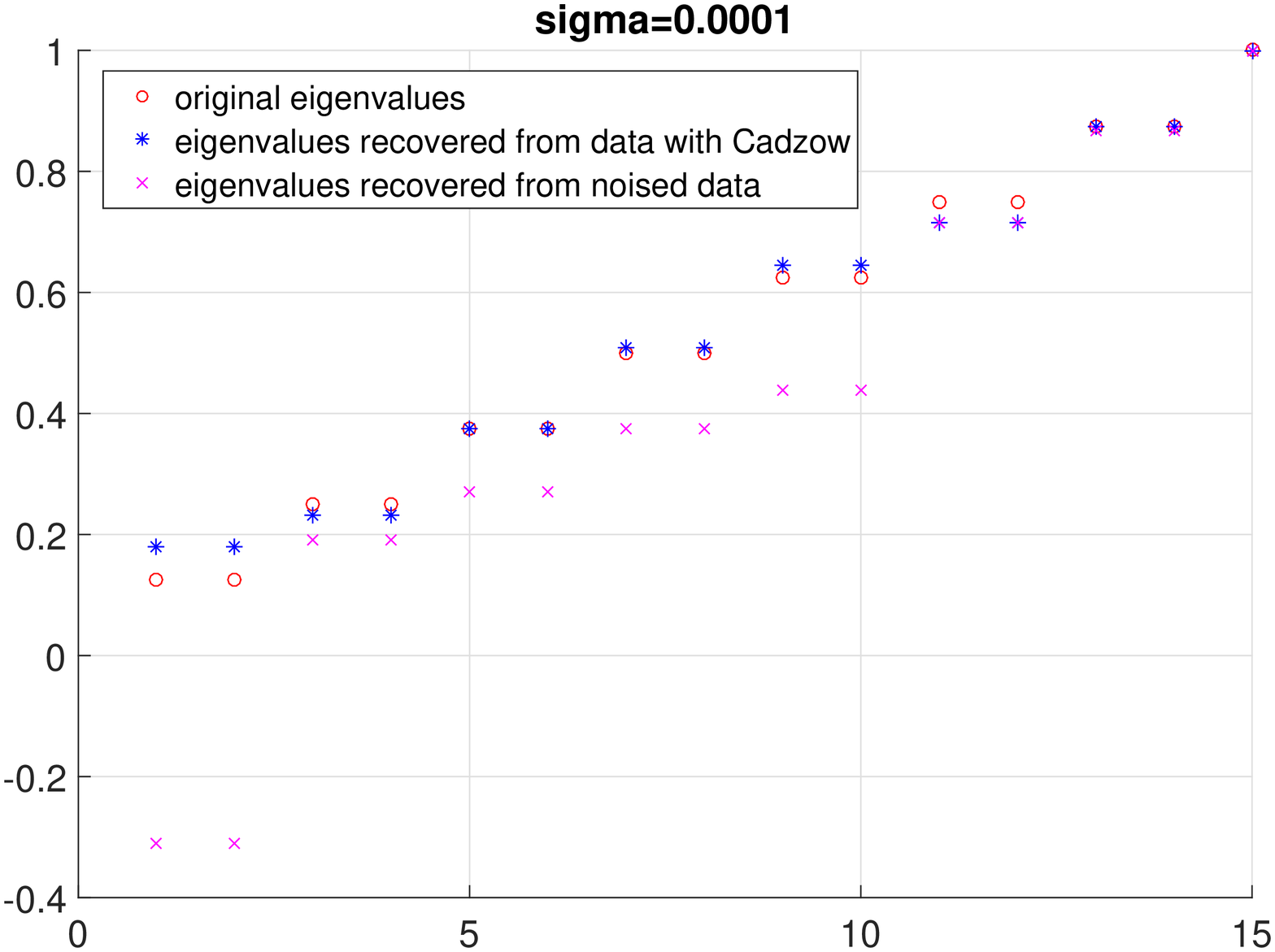}
\caption{}
\end{subfigure}
\caption{A comparison of  the spectrum reconstruction with 
and without the Cadzow denoising technique for $\sigma=10^{-4}$. }
\label{fig:4.2.1.2}
\end{figure}

\begin{figure}[tbph]
\centering
\begin{subfigure}[b]{0.49\textwidth}
\includegraphics[width=\textwidth]{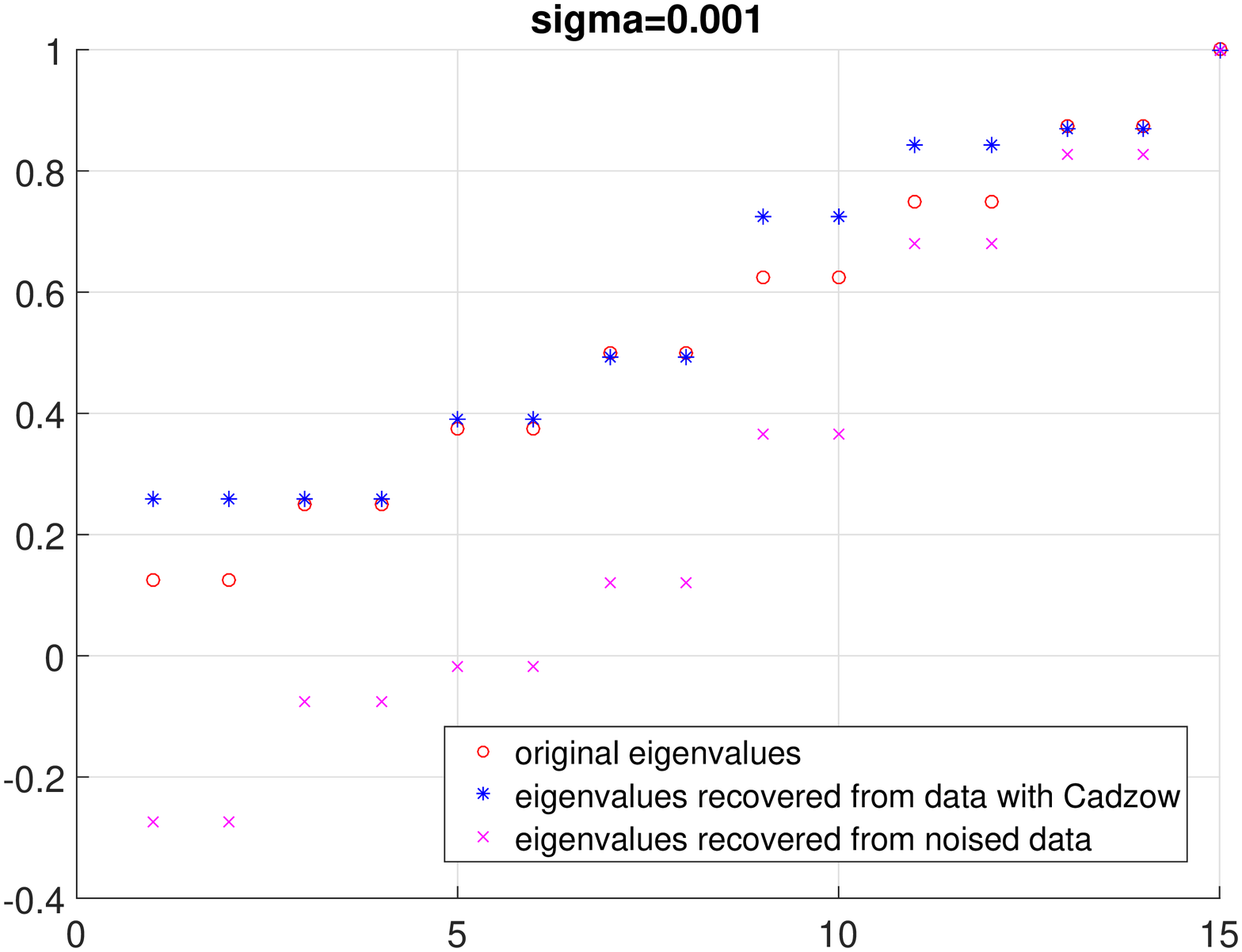}
\caption{}
\end{subfigure}
\begin{subfigure}[b]{0.49\textwidth}
\includegraphics[width=\textwidth]{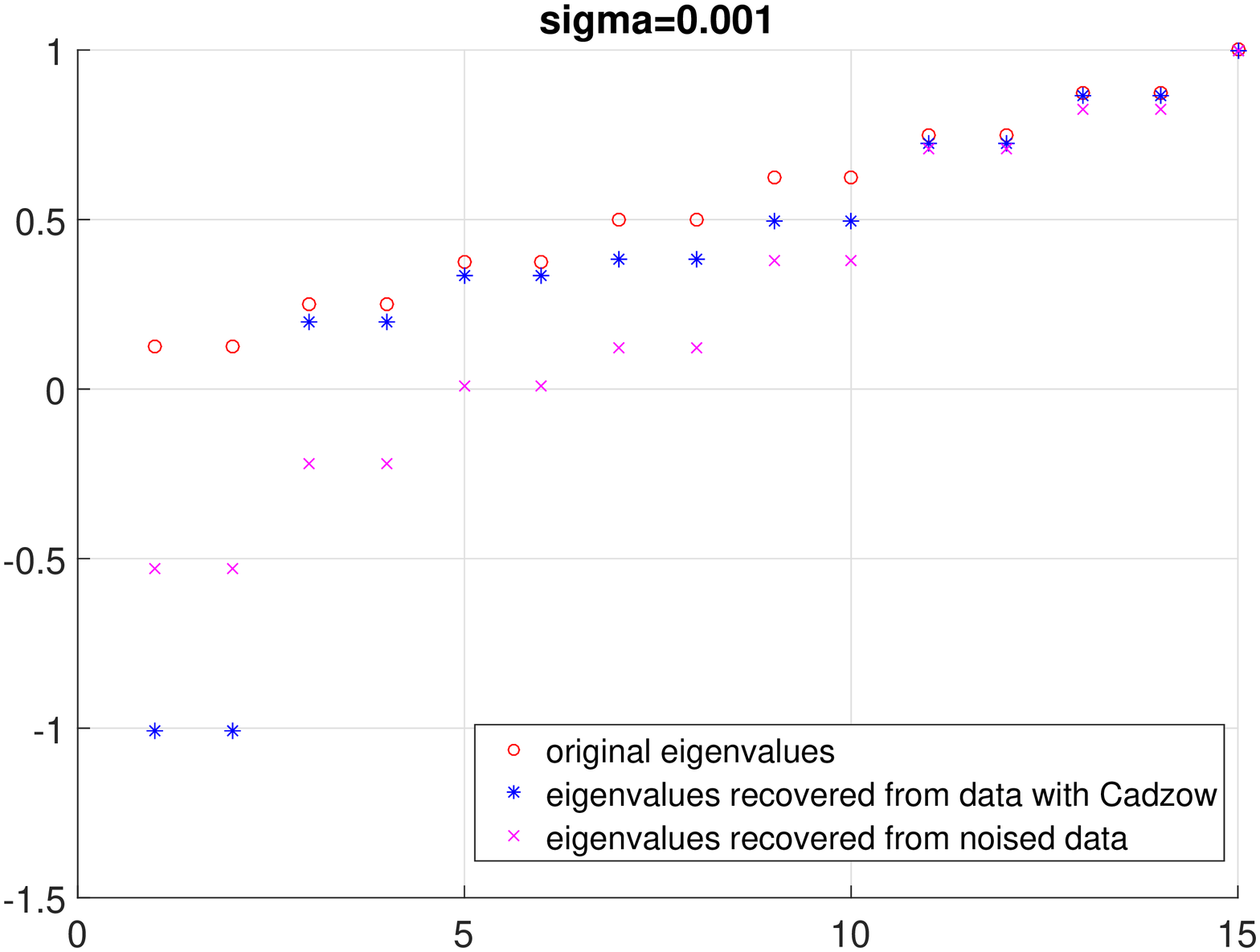}
\caption{}
\end{subfigure}
\begin{subfigure}[b]{0.49\textwidth}
\includegraphics[width=\textwidth]{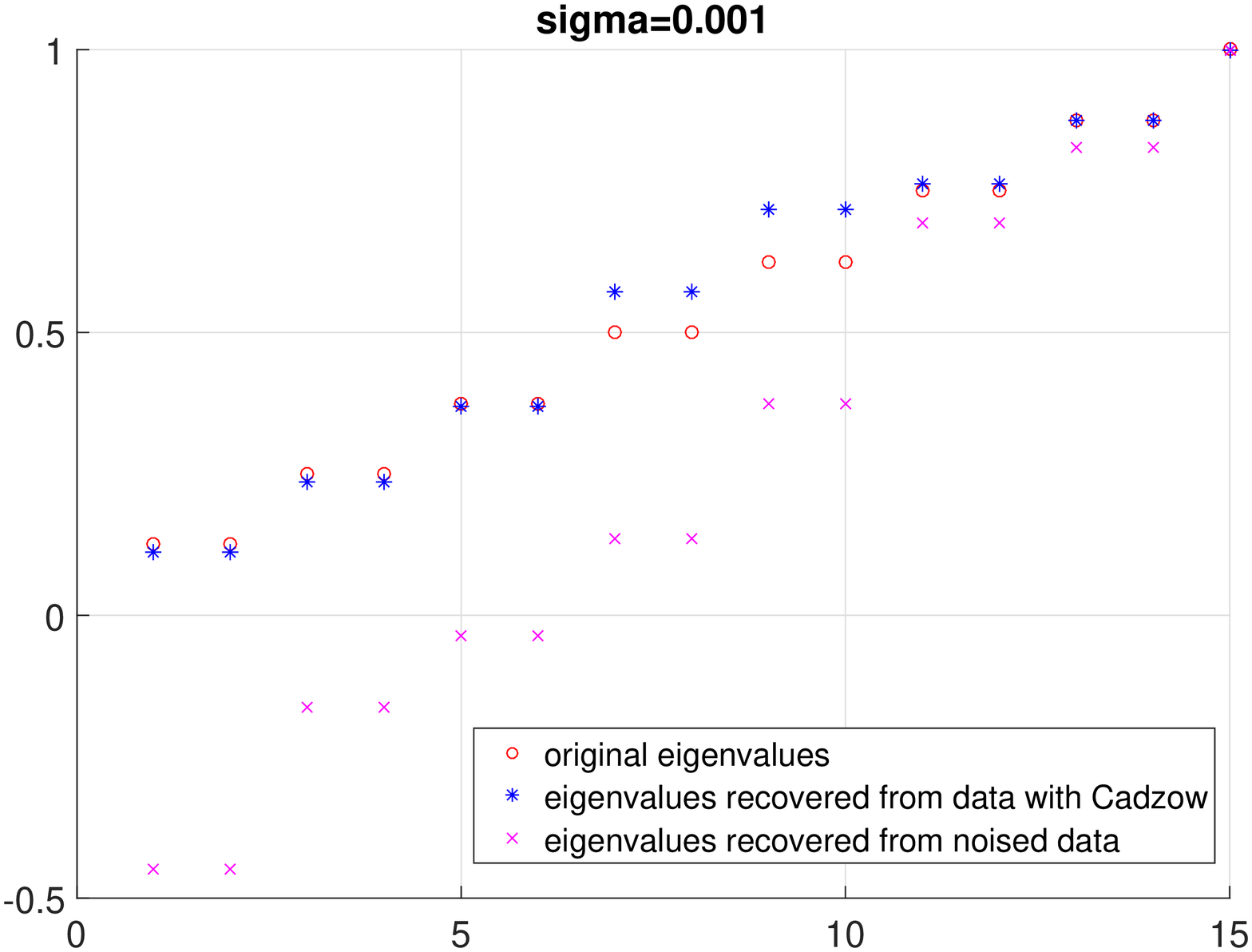}
\caption{}

\end{subfigure}
\caption{A comparison of  the spectrum reconstruction with 
and without the Cadzow denoising technique for $\sigma=10^{-3}$. }
\label{fig:4.2.1.3}
\end{figure}

In order to evaluate the impact of the Cadzow denoising technique when   reconstructing the spectrum of the convolution operator,  we conducted  a number of simulations on synthetic data. 
We repeated the same process as  in Section \ref{section sample reconstruction} until the denoised data $Z$ was generated. 
Then we used the results of Section \ref{frec} and Algorithm \ref{spectrum recovery for convolution operators} to recover the spectrum of the convolution operator  using separately  denoised data $Z$ and  noisy data  $\widetilde Y=S_m\widetilde\Pi$.  
The simulation results are shown in Figures \ref{fig:4.2.1.1}, \ref{fig:4.2.1.2}, and \ref{fig:4.2.1.3} for different noise standard deviations. Figure \ref{fig:4.2.1.1} shows the simulation results when the  standard deviation of the noise is $10^{-5}$. The curves in Figure \ref{fig:4.2.1.1} are simulation results for three different random choices of  noise. For Figures \ref{fig:4.2.1.2} and \ref{fig:4.2.1.3}, the noise has standard deviations $10^{-4}$ and $10^{-3}$, respectively. As shown in Figures \ref{fig:4.2.1.1}, \ref{fig:4.2.1.2}, and \ref{fig:4.2.1.3}, the Cadzow denoising technique can make a big difference for the spectrum recovery. 

\begin{figure}[!tbph]
\centering
\begin{subfigure}[b]{0.49\textwidth}
\includegraphics[width=\textwidth]{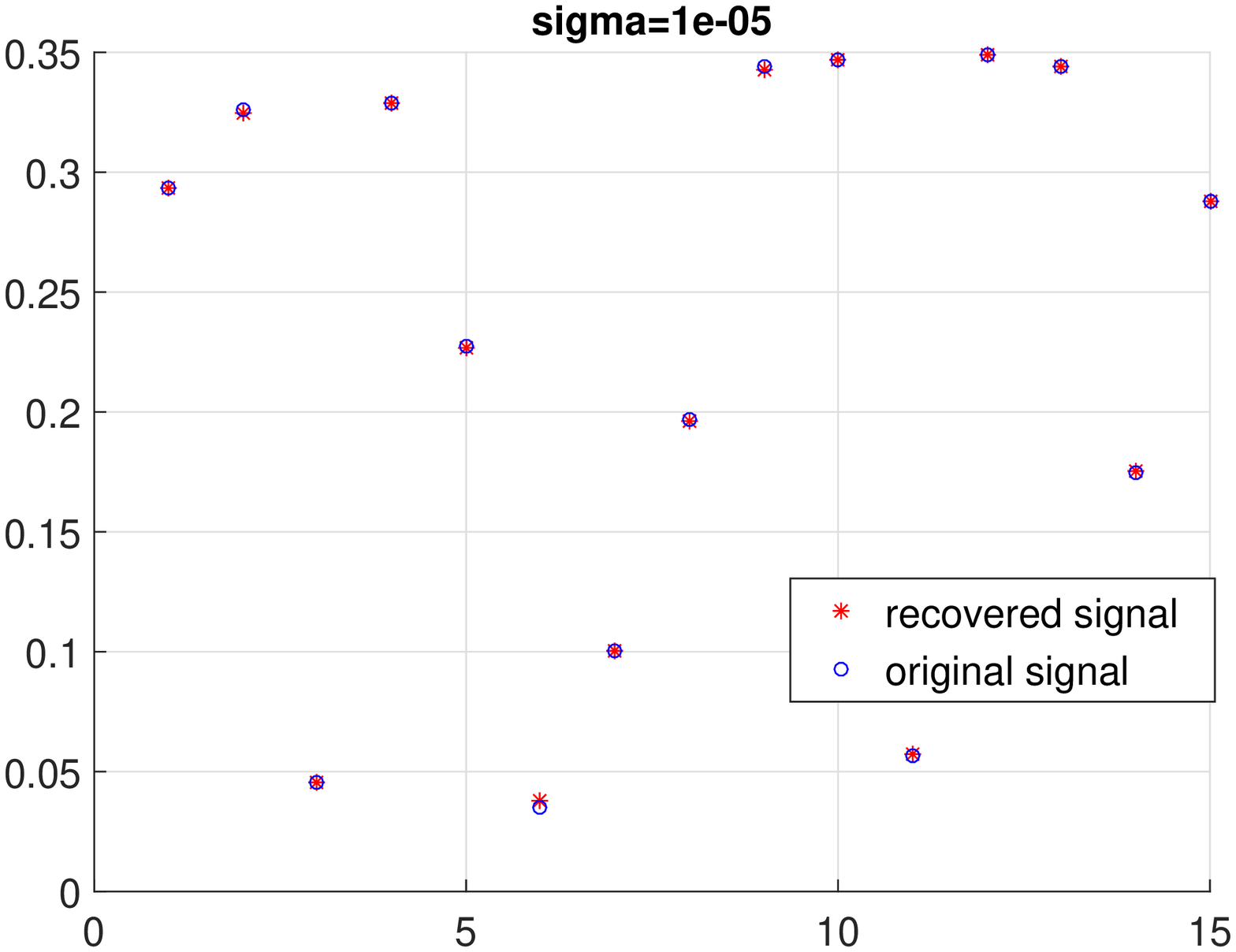}
\caption{}
\end{subfigure}
\begin{subfigure}[b]{0.49\textwidth}
\includegraphics[width=\textwidth]{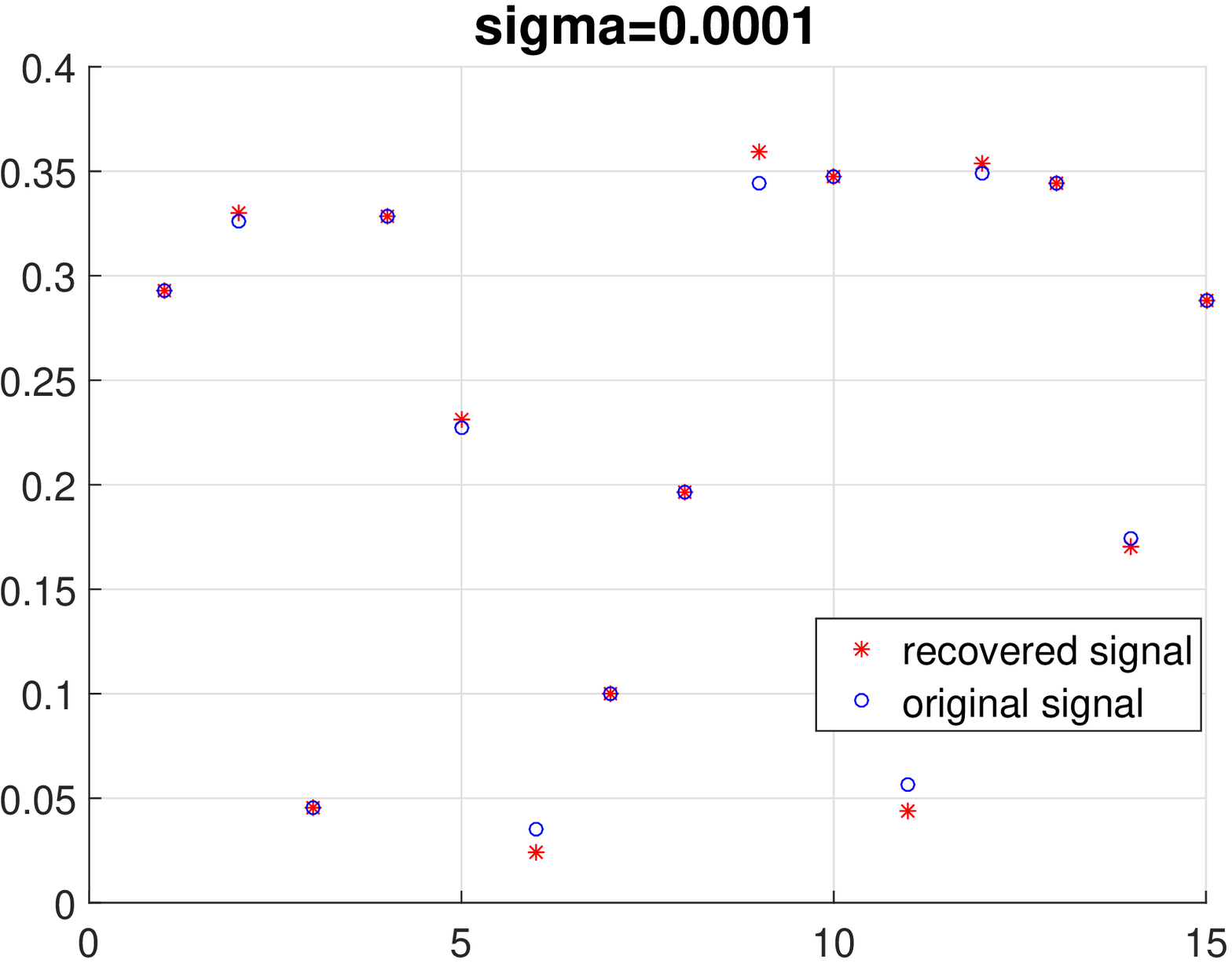}
\caption{}
\end{subfigure}
\begin{subfigure}[b]{0.49\textwidth}
\includegraphics[width=\textwidth]{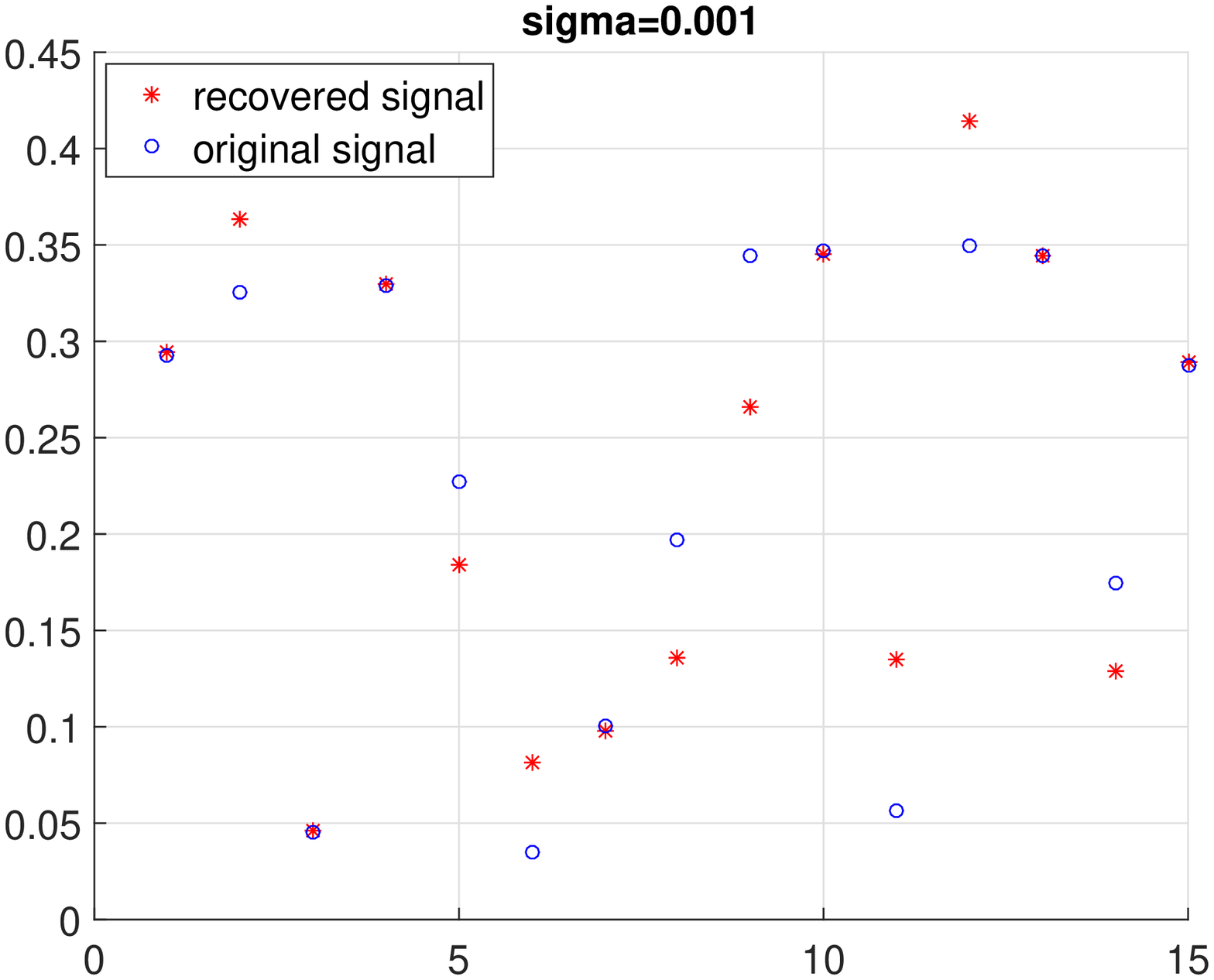}
\caption{}

\end{subfigure}
\caption{A comparison of the recovered signal and the original signal by using the estimated recovered convolution operator. }
\label{fig:4.2.1.4}
\end{figure}

Using the estimated convolution operator and the denoised data, we also evaluated the effectiveness  of  the reconstruction algorithm, i.e., Algorithm \ref{algorithm iterative}, for which the simulation  results are shown in  Figure \ref{fig:4.2.1.4}. The figure shows that if the noise is small, the recovered signals are extremely close to the original signals, which also verifies the effectiveness of the Cadzow denoising technique for  dynamical sampling. 


\subsection {Real data}\label{rd}

In this section, we describe numerical tests that we performed using   two sets of  real data. One data set documents a cooling process with a single heat source,  and the other -- a similar process with two heat sources. These data sets were labeled as ``one hotspot" and ``two hotspots", respectively. 

\begin{figure}
\centering
\includegraphics[width=0.5\textwidth]{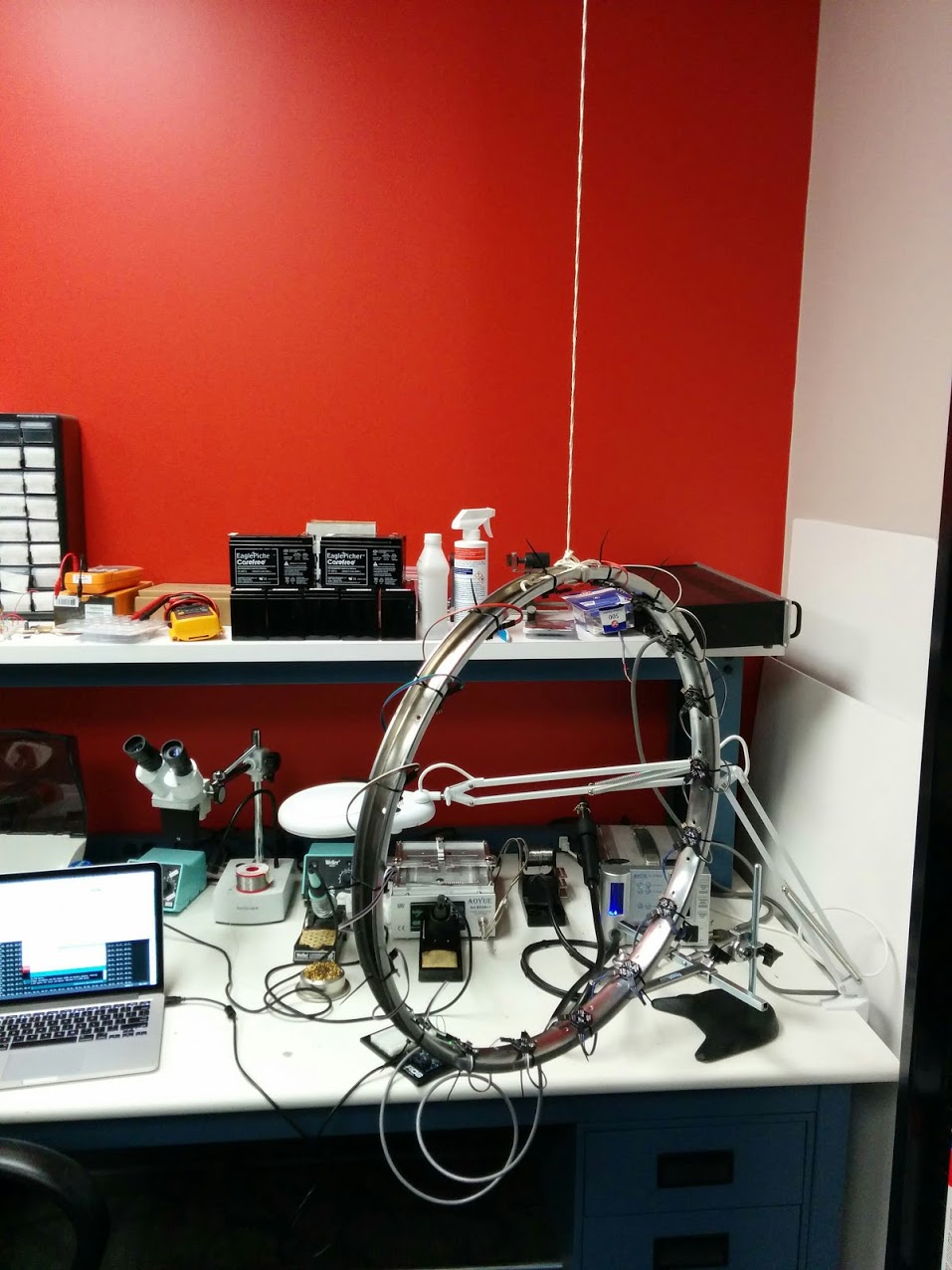}
\caption{Set-up. }\label{setup}
\end{figure}

The set-up for the real data sets is shown in Figure \ref{setup}. We used the bicycle (aluminum) wheel for the circular pattern.  Fifteen (15)  sensors are equidistantly placed around the perimeter of the wheel with 4.5 inches apart. The specified accuracy of the sensors is $0.5^\circ C$ and the temperature samples are taken  at 1.05Hz.    

The   goal was to estimate the dynamical operator and the original signals by using   information from a subset of the thermometer measuring devices, while the totality of the measurements from all devices was used as control to assess the performance of our estimations. In our reconstructions,  we did not use any a priori knowledge about the conducting material, its parameters, or  the underlying operator driving the evolution of the temperature.  Only raw, time-space subsamples of the temperatures was  used to estimate the evolution operator, and  the initial temperature distribution.  
The operator 
was assumed to be real, symmetric convolution operator whose Fourier transform consists of two monotonic pieces, so that  recovery of the spectrum of the driving operator sufficed to recover the filter.

In the experiment, the signal at time level 20 was  set as the original state. First, we smoothed the data by averaging over time to obtain a new data set $\Gamma=(\gamma_1\; \gamma_2 \ldots)$, where $\gamma_1=\sum_{i=1}^{10}f_i$, $\gamma_2=\sum_{i=11}^{20}f_i$, etc. Next, we extracted the information from the new data set at uniform locations $\Omega$ with gap $m=3$  generating the data set $S_m(\Gamma)$.  Cadzow Algorithm \ref {Cadzow denoising algorithm} is then used on $\widetilde Y=S_m (\Gamma)$ with the threshold rank  close to 2 or 3 to obtain the denoised data $Z$. Using the data $Z$, Algorithm \ref{spectrum recovery for convolution operators} was applied to estimate the  filter.  Finally, using the recovered filter,  the original signals were estimated by repeating the computations as in Section \ref{sim1}.

\begin{figure}[!tbph]
\centering
\begin{subfigure}[b]{0.49\textwidth}
\includegraphics[width=\textwidth]{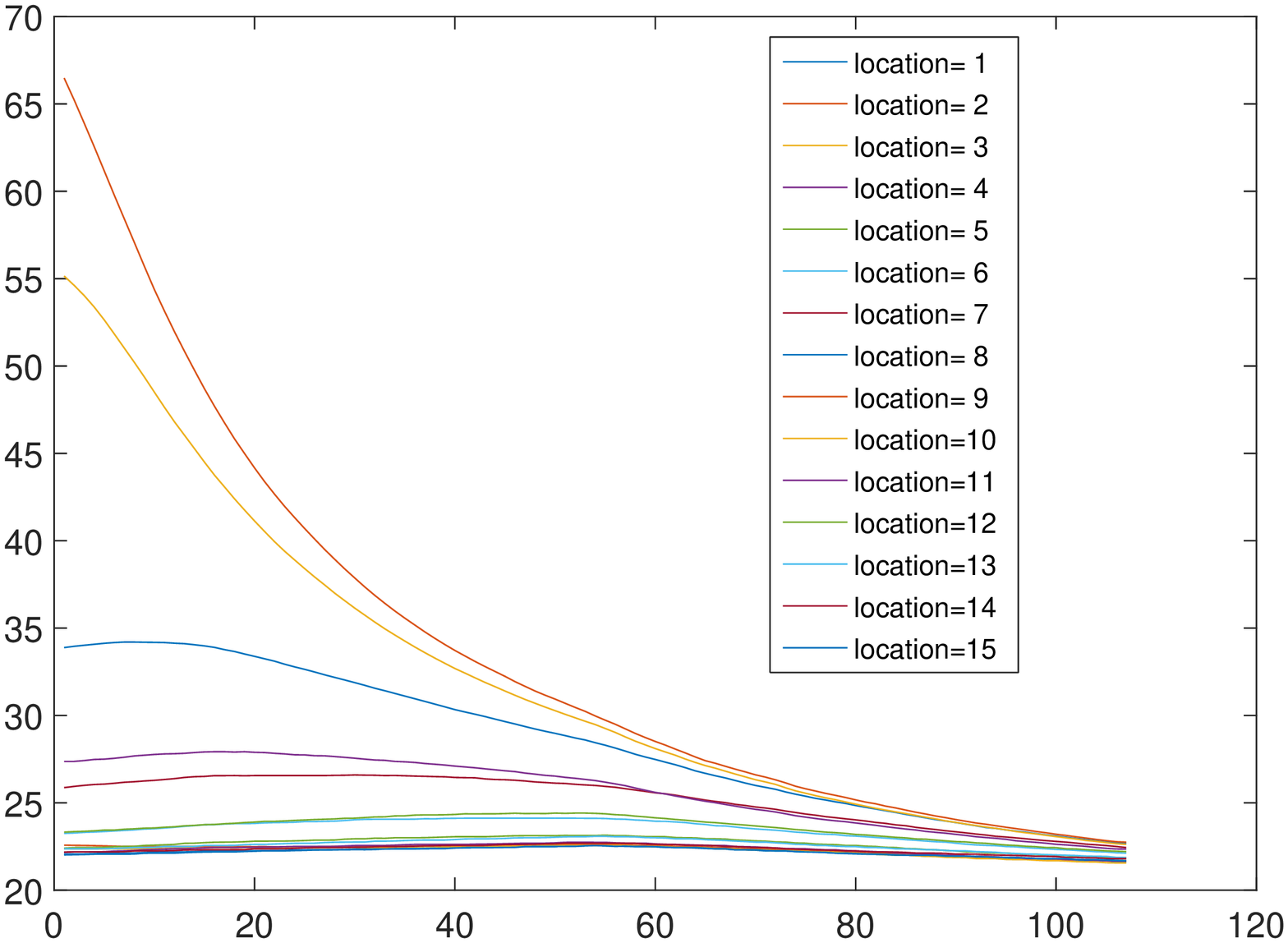}
\caption{}
\label{fig:5.1}
\end{subfigure}
\begin{subfigure}[b]{0.49\textwidth}
\includegraphics[width=\textwidth]{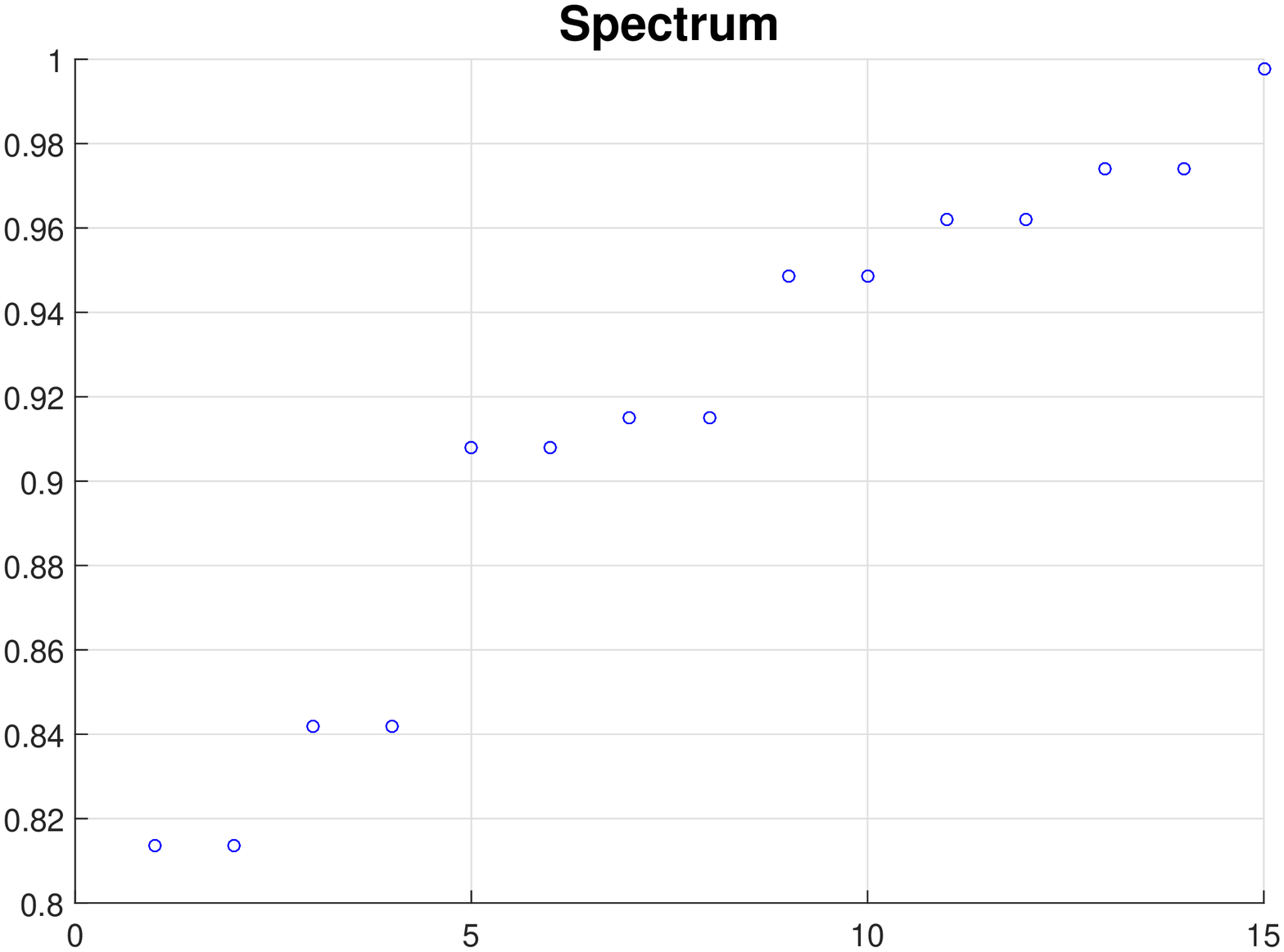}
\caption{}
\label{fig:5.2}
\end{subfigure}
\begin{subfigure}[b]{0.49\textwidth}
\includegraphics[width=\textwidth]{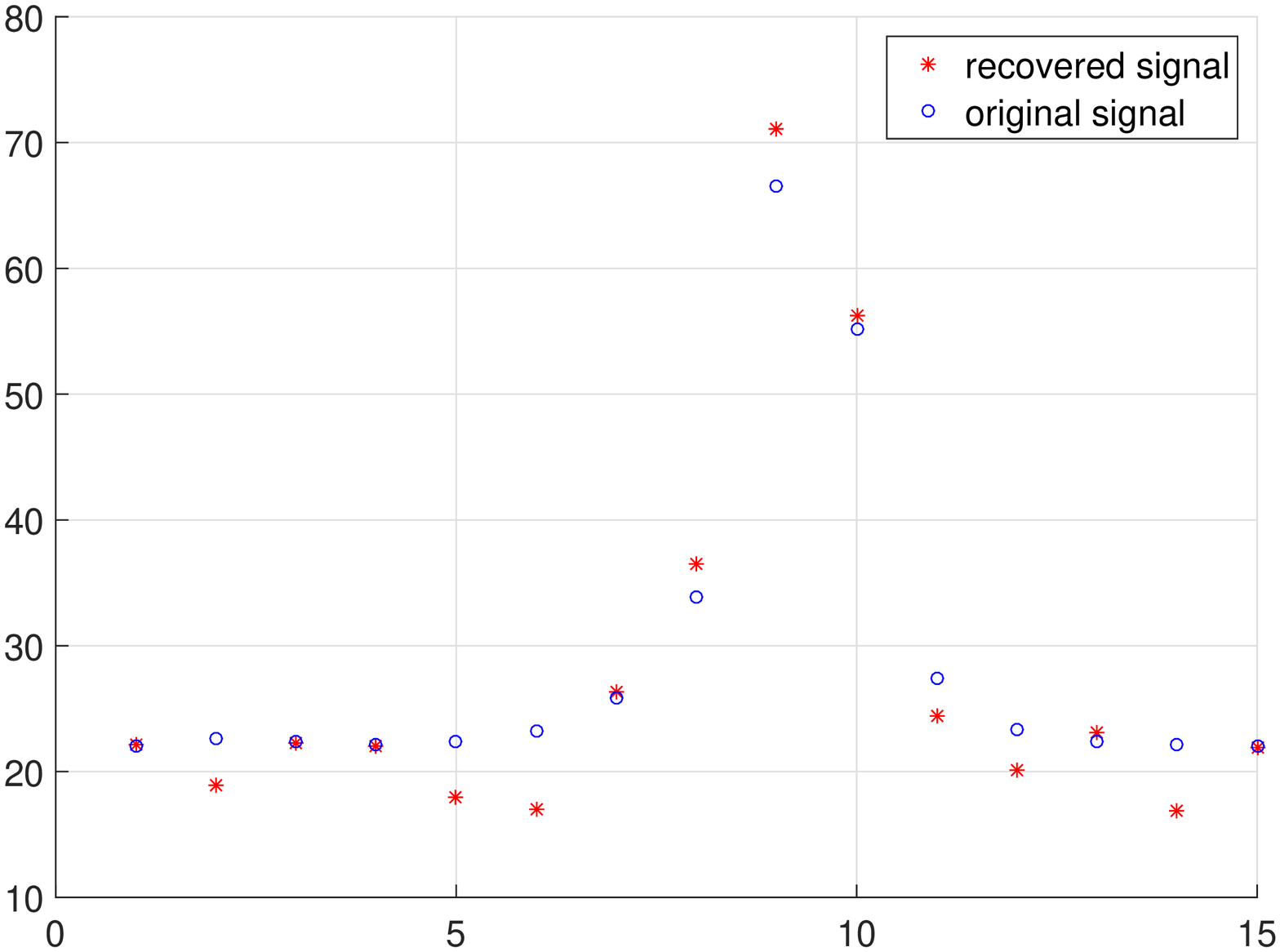}
\caption{}
\label{fig:5.3}
\end{subfigure}
\caption{Simulation results for the data set with one hotspot. Here,  (\ref{fig:5.1}) plots the evolved signals, (\ref{fig:5.2}) shows the recovered spectrum by using the data  from partial locations, and (\ref{fig:5.3}) sketches the recovered signal by using the recovered operator from partial locations and the sampled original signal.  The partial locations for recovering the operator are $\Omega=\{1,4,7,10,13\}$. To recover the original signals, we use the data from locations $\Omega_e=\{1,3,4,7,10,13,15\}$. }
\label{fig:5}
\end{figure}

The test results on the data set with one hotspot   are shown in Figure \ref{fig:5}. Figure \ref{fig:5.1} depicts the evolved signals at all 15 locations. Figure \ref{fig:5.2} shows the recovered spectrum of the evolution filter using the data from locations $\Omega=\{1,4,7,10,13\}$ to estimate the filter driving the system. Using the driving operator $A$ recovered from $\Omega$ and the necessary extra sampling locations at $\{3,15\}$ needed to recover the signal
($ \Omega_e=\Omega\cup\{3,15\}$) (see \cite {ADK13}), we reconstructed an approximation $f^\sharp$ of the signal that is  displayed in Figure \ref{fig:5.3}; it has a relative error  $\frac {\| \gamma_1-  f^\sharp\|_2} {\| \gamma_1\|_2}$ of $9.94\%$ 
compared to  the actual measurements at all 15 locations  as the reference. This relative error shows that dynamical sampling also works reasonably well for a real data set.

\begin{figure}[!tbph]
\centering
\begin{subfigure}[b]{0.49\textwidth}
\includegraphics[width=\textwidth]{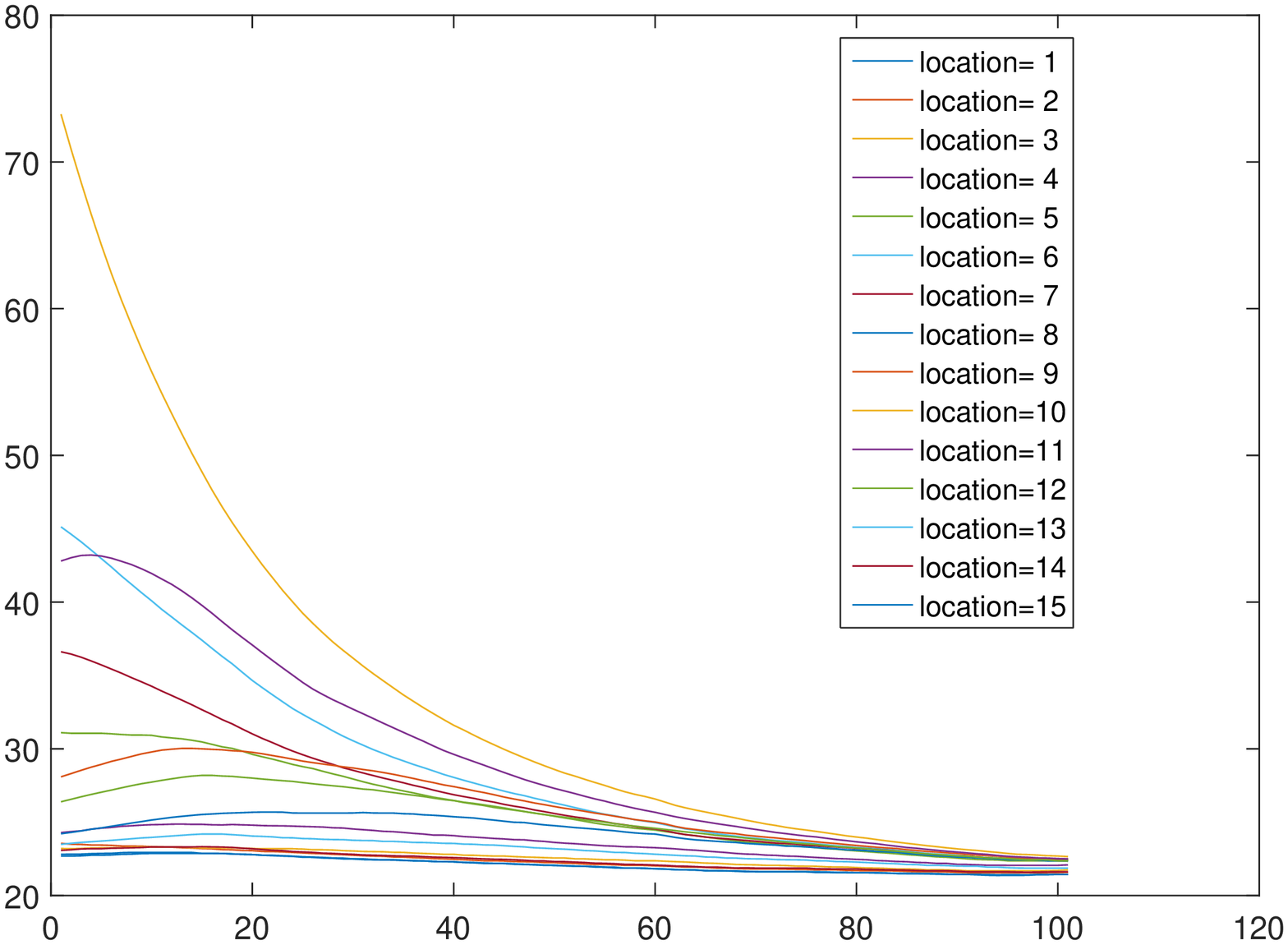}
\caption{}
\label{fig:6.1}
\end{subfigure}
\begin{subfigure}[b]{0.49\textwidth}
\includegraphics[width=\textwidth]{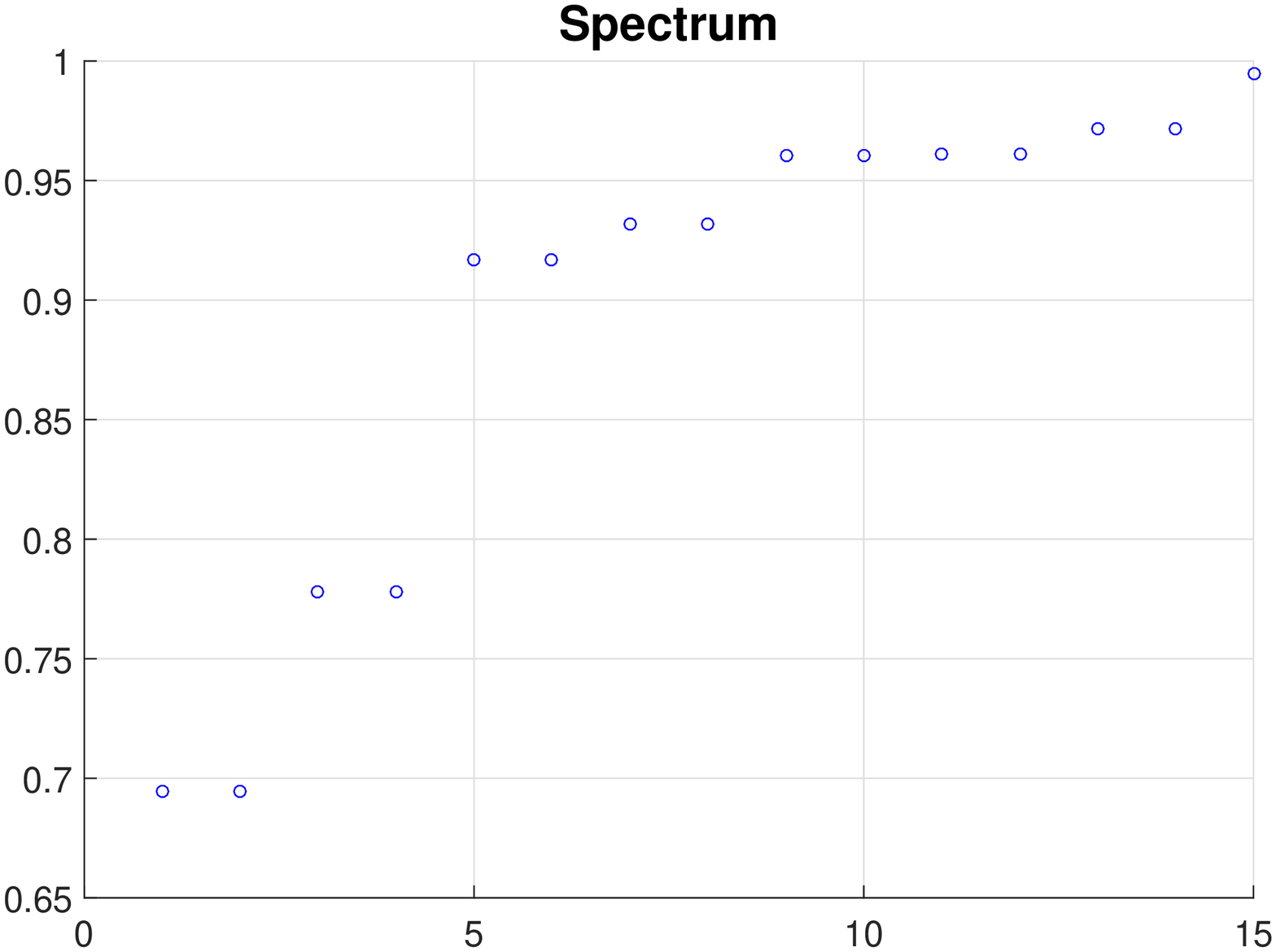}
\caption{}
\label{fig:6.2}
\end{subfigure}
\begin{subfigure}[b]{0.49\textwidth}
\includegraphics[width=\textwidth]{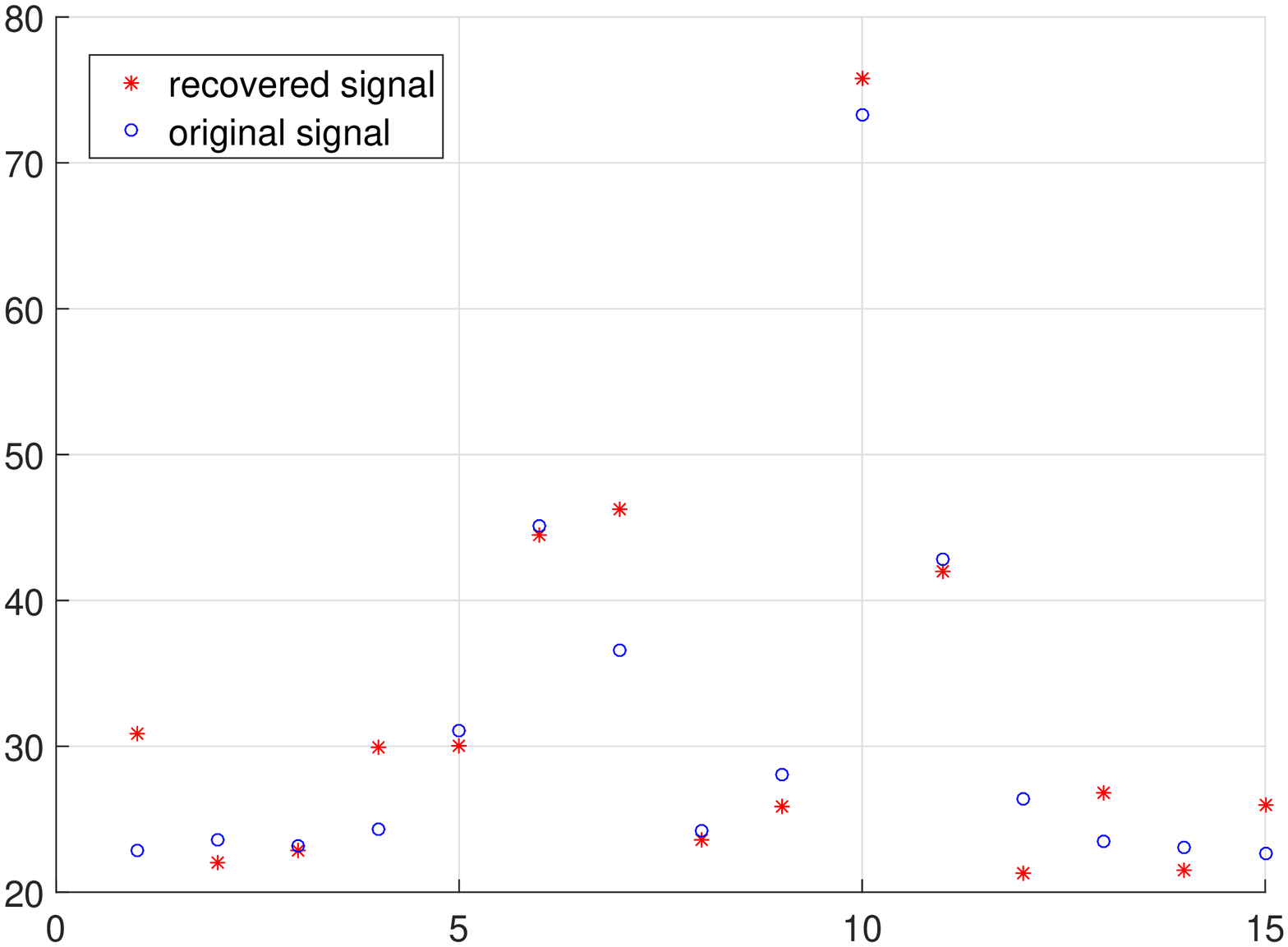}
\caption{}
\label{fig:6.3}
\end{subfigure}
\caption{Simulation results for the data set with two hotspots. Here, (\ref{fig:6.1}) plots the evolved signals, (\ref{fig:6.2}) shows the recovered spectrum by using the data from partial locations, and (\ref{fig:6.3}) sketches the recovered signal by using the recovered operator from partial locations and the sampled original signal.  The partial locations for recovering the operator are $\Omega=\{2,5,8,11,14\}$. To recover the original signals, we use the data from locations $\Omega_e=\{2,3,5,8,10,11,14\}$. }
\label{fig:6}
\end{figure}

The test results using the data set with two hotspots   are shown in Figure \ref{fig:6}. Figure \ref{fig:6.1} plots the evolved signals at the 15 locations. Figure \ref{fig:6.2} exhibits the recovered spectrum of the filter with $\Omega=\{2,5,8,11,14\}$. Using the driving operator $A$ recovered from $\Omega$ and   the  data from locations $\Omega_e=\{2,3,5,8,10,11,14\}$, we recovered an approximation of the signal that is  displayed  in Figure \ref{fig:6.3}. In this case, the relative error was $12.45\%$ compared to the actual measurements at all 15 locations. 
Such  relative error is generally considered acceptable in this kind of real applications.

\begin{figure}[!tbph]
\begin{subfigure}[b]{0.49\textwidth}
\includegraphics[width=\textwidth]{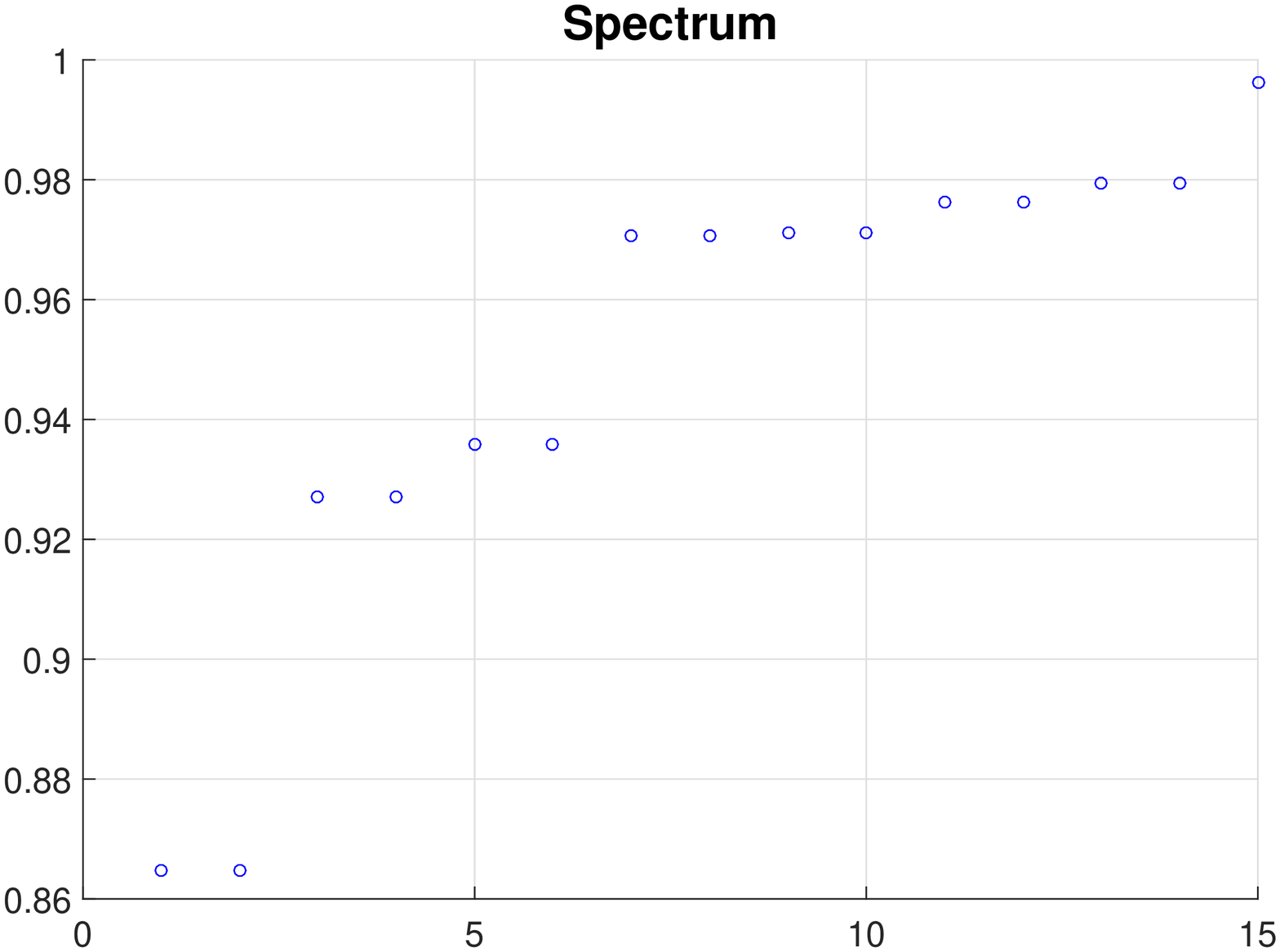}
\caption{}
\label{fig:7.1}
\end{subfigure}
\begin{subfigure}[b]{0.49\textwidth}
\includegraphics[width=\textwidth]{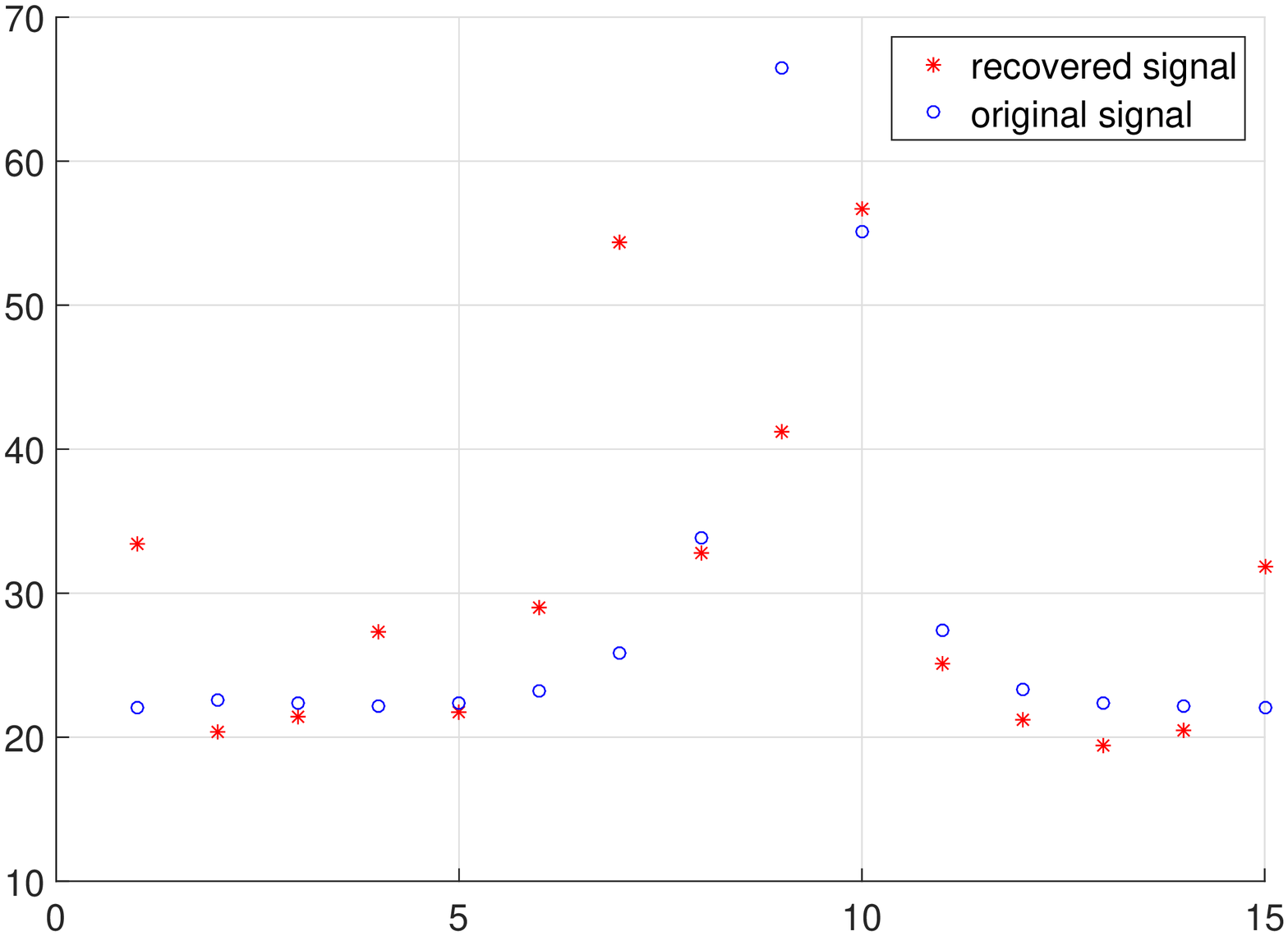}
\caption{}
\label{fig:7.2}
\end{subfigure}
\caption{Simulation results for the data set with one hotspot.  Here, (\ref{fig:7.1}) shows the recovered spectrum by using the data from  partial locations, while (\ref{fig:7.2}) plots the recovered signal by using the recovered operator from partial locations, where the partial locations for recovering the operator are $\Omega=\{2,5,8,11,14\}$. To recover the original signals, we use the samples from locations $\Omega_e=\{2,3,5,8,10,11,14\}$. }
\label{fig:7}
\end{figure}

By making similar tests with different choices of $\Omega$ and $\Omega_e$, we found that
the relative errors depend heavily on the choice of locations. The two pictures in  Figure \ref{fig:7} are the results of the same process that was used to generate the last two pictures in
Figure \ref{fig:5}. In this case, however, we chose $\Omega=\{2,5,8,11,14\}$ and $\Omega_e=\{2,3,5,8,10,11,14\}$. This choice resulted in the relative error of  $34.29\%$ which is considerably larger than the $9.94\%$ in Figure \ref{fig:5}. 

\section {Concluding remarks}
This paper introduces the problem of noise into the modeling of dynamical sampling and discusses certain unbiased linear estimators for the recovery of signals from dynamical sampling. The addition of noise to the model highlights some of the difficulties in recovering a signal from measurements in dynamical sampling, and sets the stage for more detailed studies of the information theoretic bounds and other types of estimators.

In addition, this paper studies a special case related to blind deconvolution, where the subsampling is uniform (to which  extra samples are added for the recovery of the unknown signal), and the evolution operator is unknown, but is one dimensional, symmetric, real and decreasing in the frequency domain. The existence of multiple measurements over time, along with the assumptions on the properties of the filter, allow for the recovery of the unknown signal and unknown filter; we point to some of the factors that have an adverse effect on the stability of this procedure. 

The basic algorithms and discussion of certain special cases are presented here with the intent of providing a starting point for future work on both the theoretical and algorithmic aspects of noisy instances of dynamical sampling and the case where the evolution operator is unknown. 

\vspace{13pt}
\centerline{{\bf \Large {Acknowledgement}}}
\vspace{13pt}
The research is supported by the collaborative NSF ATD grant DMS-1322099 and DMS-1322127. 
We would like to thank the organizers of SampTA 2017 in Estonia and CIMPA 2017 in Argentina for their hospitality. We would also like to thank Miklos Maroti and the anonymous reviewers for helpful comments and suggestions. Special thanks are reserved to S.~J.~Rose for imparting his wisdom, as generously as ever.

\appendix

\section {Proof of Proposition \ref{theorem2.2}}
\begin{proof}[\textbf{Proof of Proposition \ref{theorem2.2}}] It is clear that
\begin{eqnarray}
\|\epsilon_L\|_2^2&=&\left\|\left(\sum_{i=1}^{L}A_i^*A_i\right)^{-1}\left(\sum_{j=1}^{L}A_j^*\eta_j\right)\right\|^2_2\nonumber\\
&=&\sum_{j=1}^{L}\left\|\left(\sum_{i=1}^{L}A_i^*A_i\right)^{-1}A_j^*\eta_j\right\|_2^2+\sum_{j\neq k}\left\langle \left(\sum_{i=1}^{L}A_i^*A_i\right)^{-1}A_j^*\eta_j,\left(\sum_{i=1}^{L}A_i^*A_i\right)^{-1}A_k^*\eta_k\right\rangle\nonumber.
\end{eqnarray}
Since $\eta_j$ and $\eta_k$ for $j\neq k$ are independent  and mean zero, the cross terms cancel out in expectation, and one has $$E\left(\left\langle\left(\sum_{i=1}^{L}A_i^*A_i\right)^{-1}A_j^*\eta_j,\left(\sum_{i=1}^{L}A_i^*A_i\right)^{-1}A_k^*\eta_k\right\rangle\right)=0.$$
Consequently,
\begin{eqnarray}
E(\|\epsilon_L\|_2^2)=\sum_{j=1}^{L}E\left(\left\|\left(\sum_{i=1}^{L}A_i^*A_i\right)^{-1}A_j^*\eta_j\right\|_2^2\right).
\end{eqnarray}
Note that
\begin{eqnarray}
&&\left\|\left(\sum_{i=1}^{L}A_i^*A_i\right)^{-1}A_j^*\eta_j\right\|_2^2=\left\|\left(\sum_{i=1}^{L}A_i^*A_i
\right)^{-1}\sum_{l=1}^{m_j}A_j^{*(l)}\eta_j^{l}\right\|_2^2\nonumber\\
&=&\sum_{l=1}^{m_j}\left\|\left(\sum_{i=1}^{L}A_i^*A_i\right)^{-1}A_j^{*(l)}\eta_j^{l}\right\|_2^2+\sum_{l\neq p}\left\langle \left(\sum_{i=1}^{L}A_i^*A_i\right)^{-1}A_j^{*(l)}\eta_j^{l},\left(\sum_{i=1}^{L}A_i^*A_i\right)^{-1}A_j^{*(p)}\eta_j^{p}\right\rangle,\nonumber
\end{eqnarray}
where $A_j^{*(l)}$ denotes the $l$-th column of  matrix $A_j^*$ and $\eta_j^{l}$ is the $l$-th entry of $\eta_j$. 
Additionally, $\eta_j^{l}$ and $\eta_j^{p}$ are independent for $l\neq p$. It follows that
\[E\left(\left\langle \left(\sum_{i=1}^{L}A_i^*A_i\right)^{-1}A_j^{*(l)}\eta_j^{l},\left(\sum_{i=1}^{L}A_i^*A_i\right)^{-1}A_j^{*(p)}\eta_j^{p}\right\rangle\right)=0.\]
Thus, 
\begin{eqnarray}
E(\|\epsilon_L\|_2^2)&=&\sum_{j=1}^{L}\sum_{l=1}^{m_j}E\left(\left\|\left(\sum_{i=1}^{L}A_i^*A_i\right)^{-1}A_j^{*(l)}\eta_j^{l}\right\|_2^2\right)\nonumber\\
&=&\sigma^2\sum_{j=1}^{L}\sum_{l=1}^{m_j}\left\|\left(\sum_{i=1}^{L}A_i^*A_i\right)^{-1}A_j^{*(l)}\right\|_2^2\nonumber\\
&=&\sigma^2\cdot\sum_{j=1}^{L}\text{trace}\left(A_j\left(\sum_{i=1}^{L}A_i^*A_i\right)^{-2}A_j^*\right)\nonumber\\
&=&\sigma^2\cdot \text{trace}\left(\left(\sum_{j=1}^{L}A_j^*A_j\right)\left(\sum_{i=1}^{L}A_i^*A_i\right)^{-2}\right)\nonumber\\
&=&\sigma^2\cdot \text{trace}\left(\left(\sum_{i=1}^{L}A_i^*A_i\right)^{-1}\right)=\sigma^2\sum_{i=1}^{d}1/\lambda_i(L)
\end{eqnarray}
and the proposition is proved.
\end{proof}

\bibliographystyle{siam}
\bibliography{Akram_refs}

\end{document}